\documentclass{article}%
\usepackage{latexsym}
\usepackage{amsmath}
\usepackage{amsfonts}
\usepackage{amssymb}
\usepackage{graphicx}%
\usepackage{multirow}

\begin{document}

\input{amssym.def}



\baselineskip 20pt \textsf{ \textbf{{\huge \newline Penalized
contrast estimator for adaptive density deconvolution}} }

\vglue 5mm

\baselineskip 14pt

\noindent\textsf{{\large Fabienne Comte, Yves Rozenholc and Marie-Luce Taupin
}}

\vglue 5mm

\baselineskip 12pt

Density deconvolution. Model selection. Penalized contrast.
Projection Estimator.
\vglue 5mm

\noindent{\small \textit{Abstract:} The authors consider the
problem of estimating the density $g$ of independent and
identically distributed variables $X_i$, from a sample $Z_1,
\dots, Z_n$ where $Z_i=X_i+\sigma\varepsilon_i$, $i=1, \dots, n$,
$\varepsilon$ is a noise independent of $X$, with
$\sigma\varepsilon$ having known distribution. They present a
model selection procedure allowing to construct an adaptive
estimator of $g$ and to find non-asymptotic bounds for its
$\mathbb{L}_2(\mathbb{R})$-risk. The estimator achieves the minimax rate of
convergence, in most cases where lowers bounds are available.
A simulation study gives an illustration
of the good practical performances of the method.}

\vglue 5mm

\noindent\textbf{D\'econvolution adaptative de densit\'e par contraste
p\'enalis\'e.}

\vglue 5mm

\noindent{\small \textit{R\'{e}sum\'{e} :} Les auteurs
consid\`erent le probl\`eme de d\'econvolution c'est-\`a-dire de
l'estimation de la densit\'e de variables al\'eatoires
identiquement distribu\'ees $X_i$, \`a partir de l'observation de
$Z_i$ o\`u $Z_i=X_i+\sigma\varepsilon_i$, pour $i=1, \dots, n$,
o\`u les erreurs $\sigma \varepsilon_i$ sont de densit\'e connue.
Par une proc\'edure de s\'election de mod\`eles qui permet
d'obtenir des bornes de risque non asymptotiques, ils construisent
un estimateur adaptatif de la densit\'e des $X_i$. L'estimateur
atteint de fa\c{c}on automatique la vitesse minimax dans la
plupart des cas, que les erreurs ou la densit\'e \`a estimer
soient peu ou tr\`es r\'eguli\`eres.
Une \'etude par simulation illustre les bonnes performances
pratiques de la m\'ethode.}

\vglue 5mm

\bigskip\noindent\textsf{1.\ INTRODUCTION}

\bigskip\noindent We observe $Z_1, \cdots, Z_n$, $n$
independent and identically distributed (i.i.d.) copies of $Z$ in the model
$$
Z=X+\sigma\varepsilon,$$ where $X$ and $\varepsilon$ are
independent random variables, with unknown density $g$ for $X$,
known density $f_\varepsilon$ for $\varepsilon$, and known noise
level $\sigma$. In this model, we aim  at estimating the density
$g$ without any prior knowledge on its smoothness, using the
observations $Z_1,\cdots,Z_n$ and the knowledge of the convolution
kernel $\sigma f_\varepsilon(\cdot/\sigma)$. The parametrer $\sigma$ is only estimable under more restrictive conditions on
$g$, such as a lower bound on its Fourier transform. However, under the usual
conditions on $g$ (as in the currrent paper), $\sigma$ has to be known. We refer to Butucea
and Matias~(2005) for the problem of the estimation of $\sigma$ as well as for results about density deconvolution when
$\sigma$ is unknown in such a model.

In density deconvolution, two factors determine the estimation
accuracy. First, the smoothness of the density to be estimated,
$g$, and second the smoothness of the errors density, the worst
rates of convergence being obtained for the smoothest errors
density. Indeed, due to the independence of $X$ and $\varepsilon$,
the density $h$ of $Z$ is $h(\cdot)=g*(\sigma
f_\varepsilon(\cdot/\sigma))$, where $*$ denotes the convolution
product, and if $f_\varepsilon$ is very smooth then so is $h$, the
density of the observations and thus it is difficult to recover
$g$.

In this context, we consider two classes of errors: first the so called ordinary
smooth errors with polynomial decay of their  Fourier transform
and second, the supersmooth errors with Fourier transform having an
exponential decay.

Most previous results concern kernel estimators and densities $g$ to be
estimated belonging to H\"older or Sobolev classes with known order $s$.
One can cite among  others Carroll and Hall
(1988), Devroye (1989), Fan~(1991a, b), Liu and Taylor~(1989),
Masry~(1991), Stefanski and Carroll~(1990), Zhang~(1990),
Koo~(1999), Cator~(2001).

Smoother densities $g$ with exponential decrease of their Fourier
transform, have been first considered by Pensky and
Vidakovic~(1999), Butucea~(2004) and Butucea and Tsybakov~(2004).
The latter study the sharp optimality (in a minimax sense) by
using non adaptive kernel estimators and provide an adaptive
estimator in some special case. The former is the first paper
dealing with adaptivity in a general context. This first adaptive
estimator is a wavelet estimator, that achieves the minimax rates
when $g$ belongs to some Sobolev class, but that fails in reaching
the minimax rates when both the errors density and $g$ are super
smooth. Let us mention also Pensky~(2002) for the estimation of
irregular functions and Fan and Koo~(2002) who consider wavelet
estimators for densities belonging to Besov spaces. Lastly,
analogously to Hesse~(1999), Delaigle and Gijbels~(2004a,b) study
adaptive methods using cross validation and bootstrap methods in
the kernel context.

In the spirit of Barron et al. (1999), we build an adaptive
estimator $\tilde g$, constructed by model selection, and more
precisely by minimization of a penalized contrast function. We
show that $\tilde g$ is adaptive in the sense that its
construction does not require any prior smoothness knowledge on g
and that its rate of convergence is the minimax rate of
convergence (up to some logarithmic factor) in all cases where
lower bounds are previously known, that is in most cases. More
precisely, we establish non-asymptotic bounds for its integrated
quadratic risk that ensure an automatic trade-off between a bias
term and a penalty term, only depending on the observations and on
$\sigma f_\varepsilon(\cdot/\sigma)$.

The estimator automatically achieves the best rate obtained by the
collection of non-penalized estimators when the (unknown) optimal
space is selected, exactly or sometimes within a negligible
logarithmic factor. In all cases where lower bounds are available,
this best rate is the minimax rate of convergence. In particular, when both the density and the errors are super
smooth ($\delta>0$ and $r>0$ in (A$_2^{\varepsilon}$) and
(R$_1^X$) below), our adaptive estimator significantly improves the rates
given by the adaptive estimator built in Pensky and
Vidakovic~(1999) whereas both adaptive estimators have the same rate in the
other cases (see Section 4.3).

The paper is organized as follows. In Section 2, we present the
assumptions and the estimators. In Section 3 we give upper bounds
for the $\mathbb{L}_2(\mathbb{R})$-risk of the estimator, when the
smoothness of $g$ is known, and study the optimality in a minimax
sense of the resulting rates. In Section 4, we give upper bounds
of the $\mathbb{L}_2(\mathbb{R})$-risk of the penalized minimum
contrast estimator $\tilde g$ when no prior knowledge on the
smoothness of $g$ is used.
The theoretical results are illustrated by a simulation
study in Section 5, and all the proofs are gathered in Section 6.

\bigskip\noindent\textsf{2.\ CONSTRUCTION OF THE ESTIMATORS}

\bigskip\noindent For $u$ and $v$ in $\mathbb{L}_2(\mathbb{R})$,
$u^*$ denotes the Fourier transform of $u$, $u^*(x)=\int e^{itx}u(t)dt$, $u*v$
is the
convolution product, $u*v(x)=\int u(t)v(x-t)dt$, $\parallel u\parallel=\left(\int
|u|^2(x)dx\right)^{1/2}$, and $\langle s, t \rangle= \int
s(x)\overline{t(x)}dx$.

\bigskip\noindent\textit{2.1 Model and Assumptions}

\bigskip\noindent We require that $f_\varepsilon$
belongs to $\mathbb{L}_2(\mathbb{R})$ and that for all $x\in
\mathbb{R},f_\varepsilon^*(x)\not=0$. We
consider that:
$$\mbox{(A}_1^{X,\varepsilon}\mbox{):}
\mbox{ The sequences }(\varepsilon_i)_{i\in \mathbb{N}}\mbox{ and
}(X_i)_{i\in \mathbb{N}} \mbox{ are sequences of  independent random variables.}
$$
The smoothness of $f_\varepsilon$ is described by the following assumption.
\begin{eqnarray*} {\rm (A}_2^{\varepsilon}\mbox{):} &&
\mbox{ There exist nonnegative numbers } \kappa_0, \kappa_0', \gamma,~ \mu, \mbox{ and }\delta \mbox{ such
that $f_{\varepsilon}^*$ satisfies }\\ &&\kappa_0(x^2+1)^{-\gamma/ 2}\exp\{-\mu\vert x\vert^\delta\}\leq
|f_\varepsilon^*(x)|\leq \kappa_0'(x^2+1)^{-\gamma/ 2}\exp\{-\mu\vert
x\vert^\delta\}
\end{eqnarray*}

Only the left-hand side of (A$_2^{\varepsilon}$) is required for
upper bounds whereas the right-hand side is useful when we
consider lower bounds and optimality, in a minimax sense, of our
estimators.

When $\delta=0$ in (A$_2^{\varepsilon}$), the errors are usually called ``ordinary
smooth'' errors, and they are called ``super smooth'' when $\mu>0$ and $\delta>0$.
Indeed densities satisfying (A$_2^{\varepsilon}$) with $\delta>0$ and
$\mu>0$ are infinitely differentiable. The standard examples for
super smooth densities are the following: Gaussian or Cauchy
distributions are super smooth of order $\gamma=0, \delta=2$ and
$\gamma=0, \delta=1$ respectively. For ordinary smooth densities,
one can cite for instance the double exponential (also called
Laplace) distribution with $\delta=0=\mu$ and $\gamma=2$. Although
densities with $\delta>2$ exist, they are difficult to express in
a closed form. Nevertheless, our results hold for such densities.
Furthermore, the square integrability of $f_{\varepsilon}$ and
(A$_2^{\varepsilon}$) require that $\gamma> 1/2$ when $\delta=0$.

By convention, we set $\mu=0$ when $\delta=0$ and we assume that
$\mu>0$ when $\delta> 0$. In the same way, if $\sigma=0$, the
$X_i$'s are directly observed without noise and we set
$\mu=\gamma=\delta =0$ in this case.

Although, slower rates of convergence for estimating $g$ are
obtained for smoother error density, those rates can be improved
by some additional regularity conditions on $g$. Those
regularity conditions are described  as follows.
\begin{eqnarray*}
{\rm (R}_1^X{\rm ):} & & \mbox{There exists some positive real numbers } s, r, b \mbox{ such that
} g \mbox{ belongs to }\\&&
 \mathcal{ S}_{s,r,b}(C_1)=\{\psi\; /\;  \int_{-\infty}^{+\infty}
|\psi^*(x)|^2(x^2+1)^{s}\exp\{2b |x|^{r}\} dx\leq C_1\}
\\
{\rm (R}_2^X{\rm ):} & &
\mbox{ There exists } d >0 \mbox{
  such that }
\forall x \in \mathbb{R}, \;  \vert g^*(x)\vert \leq \mbox{1}\!\mbox{I}_{[-d,d]}(x).
\end{eqnarray*}
The smoothness classes described by (R$_1^X$) are classically
considered both in deconvolution and in ``direct'' density
estimation, with $\mathcal{ S}_{s,0,b}(C_1)$ known as Sobolev
classes. The densities satisfying (R$_1^X$) with $r>0, b>0 $ are
infinitely many times differentiable, admit analytic continuation
on a finite width strip when $r=1$ and on the whole complex plane
if $r=2$. The densities satisfying {\rm (R}$_2^X${\rm )}, often
called entire functions, admit analytic continuation in the whole
complex plane (see Ibragimov and Hasminskii ~(1983)).

Subsequently, the density $g$ is supposed to satisfy the following
assumption.
$${\rm (A}_3^X{\rm ):} \mbox{ The density } g\in \mathbb{L}_2(\mathbb{R}) \mbox{
  and there exists }  M_2>0, \mbox{ such that }
\int x^2g^2(x)dx <M_2<+\infty. $$
Assumption {\rm (A}$_3^X${\rm )} which is due to the
construction of the estimator,  is quite unusual in density
estimation. Nevertheless it already appears in density
deconvolution in a slightly different way in  Pensky and
Vidakovic~(1999) who assume, instead of {\rm (A}$_3^X${\rm )}  that $\sup_{x
\in \mathbb{R}} \vert x\vert g(x)<\infty$. It is important to note
that Assumption {\rm (A}$_3^X${\rm )}  is very unrestrictive.

All densities having tails of order $|x|^{-(m+1)}$ as $x$ tends to
infinity satisfy {\rm (A}$_3^X${\rm )}  only if $m>1/2$. One can cite for
instance the Cauchy distribution or all stable distributions with
exponent $r>1/2$ (see Devroye~(1986)). But, the Levy
distribution, with exponent $r=1/2$ does not
satisfies {\rm (A}$_3^X${\rm )}.

\bigskip\noindent\textit{2.2 The projection spaces }

\bigskip\noindent Consider $\varphi(x)=\sin(\pi x)/(\pi x),$ and let
$\varphi_{m,j}(x) = \sqrt{L_m} \varphi(L_mx-j)$, $m\in
\mathcal{M}_n=\{1,\cdots,m_n\}$.
It is well known (see for instance Meyer (1990), p.22) that $\{\varphi_{m,j}\}_{j \in
\mathbb{Z}}$ is an orthonormal basis of the space of square
integrable functions having a Fourier transform with compact
support included into $[-\pi L_m, \pi L_m]$. We denote by $S_m$
such a space and by
 $(S_m)_{m\in \mathcal M_n}$
 this collection of linear spaces. In other words
$$S_m= \{\sum_{j\in {\mathbb Z}} a_{m,j}\varphi_{_{m,j}}, \;
a_{m,j}\in {\mathbb R}\}=\{f\in
\mathbb{L}_2(\mathbb{R}), \mbox{ with } \mbox{supp}(f^*) \mbox{
included into }[-L_m\pi,L_m\pi ]\}.$$ When $L_m=2^m$, the basis
$\{\varphi_{m,j}\}$ is known as the Shannon basis, but we consider
here that $L_m=m$.

In this context, since $g_m=\sum_{j\in {\mathbb Z}} a_{m,j}
\varphi_{m,j}$ with $a_{m,j} = <g,\varphi_{m,j}>$, the orthogonal
projection of $g$ on $S_m$, involves infinite sums, we also
consider the truncated spaces $S_m^{(n)}$ defined as
$$S_m^{(n)}= \left\{\sum_{|j|\leq K_n} a_{m,j}\varphi_{m,j}, \; a_{m,j}\in {\mathbb R} \right\} \mbox{ where }
K_n \mbox{ is an integer.}$$ It is easy to see that,
$\{\varphi_{m,j}\}_{ \vert j\vert\leq K_n}$ is an orthonormal
basis of   $S_m^{(n)}$ and the orthogonal projection $g_m^{(n)}$
of $g$ on $S_m^{(n)}$ is given by $g_m^{(n)}=\sum_{\vert
j\vert\leq  K_n}a_{m,j}\varphi_{m,j}$ with
$a_{m,j}=<g,\varphi_{m,j}>$.

Associate this collection of models to the following contrast
function, for $t$ belonging to $S_m^{(n)}$
$$\gamma_n(t)=\frac{1}{n}\sum_{i=1}^n \left[\|t\|^2
-2u_t^*(Z_i)\right], \;\;\; \mbox{ with } \;\;\;u_t(x) = \frac 1{2
\pi} \left(\frac{t^*(-x)}{f_{\varepsilon}^*(\sigma x)}\right).$$
By using Parseval and inverse Fourier formulas we get that
$$\mathbb{E}\left[u_t^*(Z_i)\right]
=\frac 1{2\pi}\langle
\left(\frac{t^*(-.)}{f_{\varepsilon}^*(\sigma .)}\right)^*,
g*f_{\varepsilon}\rangle = \frac 1{2\pi}\langle
\frac{t^*(.)}{f_{\varepsilon}^*(-\sigma .)}, g^*
f_{\varepsilon}^*(\sigma .)\rangle = \frac 1{2\pi} \langle t^*,
g^*\rangle =\langle t, g\rangle,$$ and hence
$\mathbb{E}(\gamma_n(t))=\|t-g\|^2 -\|g\|^2 $ which is minimal when $t=g$. This shows
that $\gamma_n(t)$ suits well for the estimation of $g$.

\bigskip\noindent\textit{2.3 Construction of the minimum contrast estimators}

\bigskip\noindent Associated to the collection of models, the collection of
 non-penalized estimators $\hat g_m^{(n)}$ of $g$ is defined by
\begin{equation}\label{tronque}  \hat g_m^{(n)} = \arg\min_{t\in S_m^{(n)}} \gamma_n(t).\end{equation} By using that,
$t\mapsto u_t$ is linear, and that $\{\varphi_{m,j}\}_{ \vert j\vert
\leq K_n }$ is an orthonormal basis of $S_m^{(n)}$, we have $\hat
g_m^{(n)} = \sum_{\vert j\vert\leq K_n} \hat a_{m,j} \varphi_{m,j}
\;\; \mbox{ where }\;\; \hat a_{m,j}= n^{-1} \sum_{i=1}^n
u_{\varphi_{m,j}}^*(Z_i)$ and $\mathbb{E}(\hat a_{m,j})=
<g,\varphi_{m,j}>=a_{m,j}.$

\bigskip\noindent\textit{2.4 Construction of the minimum penalized contrast estimator}

\bigskip\noindent We aim at finding the
best model $\hat m$ in $\mathcal{M}_n$, based on the data and not
on prior knowledge on the smoothness of $g$, such that the risk of the resulting
estimator is almost as good as the risk of the best estimator in
the family. The model selection is performed in an automatic way,
using the following penalized criteria
\begin{equation}\label{estitronc}
\tilde g=\hat g^{(n)}_{\hat m} \mbox{ with } \hat m= \arg\min_{m\in {\mathcal M}_n}
\left[\gamma_n(\hat g_m^{(n)}) + \; {\rm pen}(m)\right],
\end{equation}
where the penalty function pen$(m)$ is defined by
\begin{equation}\label{penalite}
{\rm pen}(m)= 2a(\lambda_1+\mu \sigma^{\delta}\pi^{\delta}\lambda_2 )\frac{L_m^{\max ( 0, \min(3\delta/2-1/2,\delta) )} \Gamma(m)}n,
\end{equation}
the constant $a$ is a fixed universal constant (to be found by simulation experiments),
\begin{eqnarray}
\label{lambda1}
\lambda_1(\gamma,\kappa_0,\mu,\sigma,\delta)=\frac{(\sigma^2\pi^2+1)^{\gamma}}{\pi^{\delta}
\kappa_0^{2}R(\mu,\delta,\sigma)},
R(\mu,\delta, \sigma)= \mbox{1}\!\mbox{I}_{\{\delta=0\}} + 2\mu\delta\sigma^{\delta} \mbox{1}\!\mbox{I}_{\{0<\delta \leq 1\}} +
2\mu\sigma^{\delta}\mbox{1}\!\mbox{I}_{\{\delta>1\}},\end{eqnarray}
$$\label{lambda2}
\lambda_2=\lambda_1^{1/2}(1+\sigma^2\pi^2)^{\gamma/2}\|f_{\varepsilon}\|\kappa_0^{-1}(2\pi)^{-1/2}
\mbox{1}\!\mbox{I}_{\{1/3\leq \delta\leq 1\}} + \lambda_1
\mbox{1}\!\mbox{I}_{\{\delta >1\}},$$
\begin{equation}\label{gammadem}
\begin{array}{lc}
\; \hspace{-3cm} \mbox{ and  }~~~~~~ &
\Gamma(m)= L_m^{(2\gamma+1-\delta)}
\exp\left\{2\mu\sigma^\delta\pi^{\delta} L_m^{\delta}\right\}.
\end{array}
\end{equation}
Since $\sigma$ and $f_\varepsilon$  are known, the constants $\sigma$ and $\mu, \delta,
\kappa_0, \gamma$  defined in (A$_2^{\varepsilon}$) are also known.

\bigskip\noindent\textsf{3.\ RATES OF CONVERGENCE OF THE MINIMUM CONTRAST
ESTIMATORS $\hat g_m^{(n)}$}

\bigskip\noindent\textit{3.1 Bias-variance decomposition of
risk of $\hat g_m^{(n)}$}

\bigskip\noindent Let us first study the rate of
convergence of one estimator $\hat g_{m}^{(n)}$, when the
smoothness of $g$ is known.\newline

\noindent\textsc{Proposition 1}. \textit{ Under Assumption {\rm
(A}$_3^X${\rm )}, denote by $\Delta_1(m)=L_m\int_{-\pi}^{\pi}
\left\vert {f_\varepsilon^*(L_m
x\sigma)}\right\vert^{-2}dx/(2\pi).$ Then $ {\mathbb E}\|g-\hat
g_m^{(n)}\|^2\leq
\|g-g_m\|^2+(\pi L_m)^2(M_2+1)/K_n+2\Delta_1(m)/n.$}\newline

\noindent\textsc{Remark 1}. We point out that the
$\{\varphi_{m,j}\}$ are ${\mathbb R}$-supported (and not compactly
supported) so that we obtain an estimation on ${\mathbb R}$ and
not only on a compact set as for usual projection estimators. This
is a great advantage of this basis. Nevertheless it induces the
residual term $(\pi L_m)^2(M_2+1)/K_n$, due to the truncation $\vert
j\vert \leq K_n$. But the most important thing is that the choice
of $K_n$ does not influence the other terms. Consequently, it is
easy to check that we can find a relevant choice of $K_n$
($K_n\geq n$ under {\rm (A}$_3^X${\rm )}, that makes this last
supplementary term unconditionally negligible with respect to the
others. The choice of large $K_n$ does not change the efficiency
of our estimator from a statistical point of view but only changes
some practical computations.\newline

Let us comment the three terms in the bound of the risk.
The variance term $\Delta_1(m)/n$
depends on the rate of decay of the Fourier transform of
$f_\varepsilon$, with larger variance for smoother
$f_\varepsilon$.  Under (A$_2^{\varepsilon}$), by applying Lemma 3 in Section 6.3,
we get that
$\Delta_1(m)\leq 2\lambda_1\Gamma(m)$
where $\Gamma(m)$ is given by (\ref{gammadem}) and $\lambda_1=\lambda_1(\gamma,\kappa_0,\mu,\sigma,\delta)$ is
given by (\ref{lambda1}). In order to ensure that $\Gamma(m_n)/n$ is bounded, we only consider
$L_m=m\leq m_n$ with
\begin{eqnarray}
\label{mn}
m_n\leq \left\{
\begin{array}{ll}
\pi^{-1}n^{1/(2\gamma+1)} &\mbox{ if }\delta=0\\
\displaystyle\pi^{-1}\left[\frac{\ln(n)}{2\mu\sigma^\delta}+\frac{2\gamma+1-\delta}{2\delta\mu
  \sigma^\delta}\ln\left(\frac{\ln(n)}{2\mu\sigma^\delta}\right)\right]^{1/ \delta} &\mbox{ if
} \delta>0.
\end{array}
\right.
\end{eqnarray}
Under {\rm (A}$_3^X${\rm )} and (A$_2^{\varepsilon}$), if $K_n\geq n$, then we have
\begin{eqnarray}
{\mathbb E}\|g-\hat g_m^{(n)}\|^2&\leq&
\|g-g_m\|^2+2\lambda_1\Gamma(m)/n+(\pi L_m)^2(M_2+1)/n\label{vit2}
\end{eqnarray}

Finally, since $g_m$ is the orthogonal
projection of $g$ on $S_m$,  we get that $g_m^*=g^*\mbox{1}\!\mbox{I}_{[-L_m\pi, L_m\pi]}$
and therefore $ \|g-g_m\|^2
=(2\pi)^{-1} \|g^*-g_m^*\|^2 =
(2\pi)^{-1} \int_{|x|\geq \pi L_m} |g^*|^2(x)dx.$

\bigskip\noindent\textit{3.2 Order of the risk of $\hat g_m^{(n)}$ under regularity assumptions on $g$}

\bigskip\noindent Under {\rm (R}$_2^X${\rm )} and  {\rm (A}$_3^X${\rm )}, by choosing
$\pi L_m=d$, and $K_n\geq n$, the bias term $\parallel
g-g_m\parallel^2=0$, the bound (\ref{vit2}) becomes ${\mathbb
E}(\|g-\hat g_m^{(n)}\|^2) \leq 2\lambda_1 d^{(2\gamma+1-\delta)}
\exp\left\{2\mu\sigma^\delta\pi^{\delta}
d^{\delta}\right\}/n+d^2(M_2+1)/(\pi^2 n), $ and the density $g$
is estimated with the parametric rate of convergence.
 We refer to Ibragimov and Hasminskii~(1983) for  similar
result on the ``direct'' estimation of a density $g$ satisfying
Assumption {\rm (R}$_2^X${\rm )}, using the
observations $X_1,\cdots, X_n$.

\bigskip\noindent If now $g$ satisfies {\rm (R}$_1^X${\rm )},
$ \|g-g_m\|^2 \leq [C_1/(2\pi)](L_m^2\pi^2+1)^{-s}\exp\{-2b
\pi^{r}L_m^{r}\}$. According to (\ref{vit2}), under {\rm
(A}$_3^X${\rm )}  with $K_n\geq n$, the risk of $\hat g_m^{(n)}$
is bounded by

$$C_1(2\pi)^{-1}(L_m^2\pi^2+1)^{-s}\exp\{-2b\pi^{r}L_m^{r}\} +
2\lambda_1L_m^{(2\gamma+1-\delta)}
\exp\left\{2\mu\sigma^\delta\pi^{\delta}
L_m^{\delta}\right\}/n+(\pi L_m)^2(M_2+1)/n .$$ The optimal choices of
$L_m$ and the resulting rates are given in Table~1, for
different types of smoothness of the unknown density $g$ and
different types of known error density $f_\varepsilon$.

\begin{table}[ptbh]
\caption{Optimal choice of the length ($L_{\breve m}$)
and resulting (optimal) rates under Assumptions (A$_2^{\varepsilon}$)
and {\rm (R}$_1^X${\rm )}.}\label{rates}
\begin{center}
{\small
\begin{tabular}{clcc}\cline{3-4}\cline{3-4}
\multicolumn{2}{c}{} &\multicolumn{2}{c}{$f_\varepsilon$} \\
\multicolumn{2}{c}{} & $\delta=0$ & $ \delta>0$ \\
\multicolumn{2}{c}{} & ordinary smooth & supersmooth \\\hline
\multirow{8}{.2cm}{\\\vfill\null $g$} & $\;$ & $\;$ & $\;$ \\
& $\begin{array}{l}
  r=0\\
  \small{\mbox{Sobolev}(s)}
\end{array}$ &
$\begin{array}{l}
  \pi L_{\breve m}=O(n^{1/(2s+2\gamma +1)})\\
  \mbox{rate}=O(n^{-2s/(2s+2\gamma+1)})\\
\mbox{{\it minimax rate}}
\end{array}$  &
$\begin{array}{l}
  \pi L_{\breve m}=[\ln(n)/(2\mu\sigma^{\delta}+1)]^{1/\delta}\\
  \mbox{rate}=O( (\ln(n))^{-2s/\delta})\\
\mbox{{\it minimax rate}}
\end{array}$ \\
\cline{2-4}
& $\begin{array}{l}
  r>0\\
  \mathcal{C}^\infty
\end{array}$ &
$\begin{array}{l} \\
  \pi L_{\breve m}=\left[{\ln(n)/2b}\right]^{1/r} \\
  \mbox{ rate}= \displaystyle  O\left(\frac{\ln(n)^{(2\gamma+1)/r}}n\right)\\
\mbox{{\it minimax rate}} \\
  \;\; \end{array}$ &
$\begin{array}{c}
  L_{\breve m}  \mbox{ solution of } \\
  {L_{\breve m}}^{2s+2\gamma+1-r}\exp\{2\mu \sigma^\delta
 (\pi L_{\breve{m}})^\delta+2b \pi^r {L_{\breve m}}^r\}\\
  \qquad= O(n)\\
\mbox{{\it minimax rate if }}r<\delta \mbox{ \it and } s=0
\end{array}$\\
\hline
\end{tabular}}
\end{center}
\end{table}

Let us emphasize that the rate for $r>0, \delta>0$ is not
explicitly given, but is only written  the solution
$L_{\breve m}$ of the equation \begin{equation}\label{implicit}
{L_{\breve m}}^{2s+2\gamma+1-r}\exp\{2\mu \sigma^\delta
 (\pi L_{\breve{m}})^\delta+2b \pi^r {L_{\breve m}}^r\}= O(n).
 \end{equation}
The study of this case is of most importance since the case
$\delta>0$ contains the most studied case of Gaussian errors.
The association $\delta>0$ and $r=0$ leads usually people to conclude that this
problem is without hope when $\delta>0$ since the rates, of logarithmic order,
are indeed very slow in that case. But if we 
associate $\delta>0$ to $r>0$, then much faster than logarithmic rates are recovered (see Section 3.4). The empirical
experiments of Section 5 illustrate  that the estimation
algorithm works well in that case. Lastly, we can mention that, in the context of stochastic volatility models seen as processes observed with errors, 
most stationary distributions of standard diffusion models studied
by Comte and Genon-Catalot~(2005) happen to belong to this
class.

\bigskip\noindent\textit{3.3 About the solution of Equation
(\ref{implicit}), in the case $r>0, \delta>0$}

\bigskip\noindent The special case $r=\delta>0$ leads to the
explicit solution \begin{equation}\label{nudelopt} \pi L_{\breve
m}=[\ln(n/\ln(n)^{a})/(2\mu\sigma^\delta+2b)]^{1/r} \mbox{ with }
a=(2s+2\gamma-r+1)/{r}\end{equation} and to the rate
$[\ln(n)]^{a'}n^{-a'/(a'+\mu\sigma^\delta)}$ with $a'= (-2s \mu
\sigma^\delta + (2\gamma-r +1)b)/(r(\mu\sigma^\delta+b)).$

If $r>0$, $\delta>0$ and $r\not=\delta$, the expression of optimal
parameter $L_{\breve m}$, solution of the Equation
(\ref{implicit}), has not one single form for general $r>0$ and
$\delta>0$.

\indent \textbf{-} When $0<r <\delta$, we can precise here the
order of the rate by using some additional information on the
ratio $r/\delta<1$. We have to distinguish if $r/\delta\leq 1/2$
or $1/2<r/\delta\leq 2/3$,~\dots. More precisely, if $r/\delta\leq
1/2$, the optimal choice $L_{\breve m}$ is $$\pi L_{\breve{m}} =
\left[\frac{\ln(n)}{2\mu\sigma^\delta}
-\frac{2b}{2\mu\sigma^\delta}\left(\frac{\ln(n)}{2\mu\sigma^\delta}\right)^{r/\delta}-c
\ln\left(\frac{\ln(n)}{2\mu\sigma^\delta}\right)\right]^{1/\delta}
\mbox{ with } c= \frac{2\gamma-r+2s +1}{2\mu\sigma^\delta\delta}$$
and the rate is of order $$\ln(n)^{-2s/\delta} \exp[-2b
({\ln(n)}/(2\mu\sigma^\delta))^{r/\delta}].$$

If $1/2<r/\delta\leq 2/3$ the optimal choice of $\pi
L_{\breve{m}}$ is
$$\pi L_{\breve{m}} =
\left[\frac{\ln(n)}{2\mu\sigma^\delta}-\frac{2b}{2\mu\sigma^\delta}\left(\frac{
\ln(n)}{2\mu\sigma^\delta}\right)^{r/\delta}
+\frac{r}{\delta}\frac{(2b)^2}{2\mu\sigma^\delta}
\left(\frac{\ln(n)}{2\mu\sigma^\delta}\right)^{2r/\delta-1}
-c\ln\left(\frac{\ln(n)}{2\mu\sigma^\delta}\right)\right]^{1/\delta}$$
with the same $c$ as above, which gives the rate
$$ \ln(n)^{-2s/\delta} \exp\left[-2b
\left(\frac{\ln(n)}{2\mu\sigma^\delta}\right)^{r/\delta}+\frac{(2b)^2}{2\mu\sigma^\delta}
\frac{r}{\delta}\left(\frac{\ln(n)}{2\mu\sigma^\delta}\right)^{2r/\delta-1}\right].$$
If $2/3<r/\delta\leq 3/4$, we have another choice of $\pi
L_{\breve m}$ with another rate.

\indent \textbf{-} When $0<\delta< r$, we can also precise the
order of the rate of our estimator, by using, once again, some
additional information on the ratio $\delta/r.$ For instance, if
$\delta/r\leq 1/2$, the optimal choice $L_{\breve m}$ is
$$\pi L_{\breve{m}} =
\left[\frac{\ln(n)}{2b}
-\frac{2\mu\sigma^\delta}{2b}\left(\frac{\ln(n)}{2b}\right)^{\delta/r}-c
\ln\left(\frac{\ln(n)}{2b}\right)\right]^{1/r} \mbox{ with } c=
\frac{2\gamma-r+2s +1}{2br}$$ and the rate is of order
\begin{equation}\label{comparevit} \ln(n)^{(2\gamma+1-\delta)/r} \exp[2\mu
\sigma^{\delta}({\ln(n)}/(2b))^{\delta/r}]/n.\end{equation}
As in the case $0<r<\delta$, we obtain  a different rate for $1/2<\delta/r<2/3.$\newline

It follows that in the case $r>0$ and $\delta>0$, the rate depends
on the integer $k$ such that $r/ \delta$ or $\delta/r$ belongs to
the interval $I_k=]k/(k+1); (k+1)/(k+2)]$. We are , to our
knowledge, the first ones to have noticed this (unavoidable)
particularity of the rates.

\bigskip\noindent\textit{3.4 About the optimality of $\hat g_m^{(n)}$ when $g$
belongs to $\mathcal{ S}_{s,r,b}(C_1)$}

\bigskip\noindent The rates $n^{-2s/(2s+2\gamma +1)}$ ($\delta=0,
r=0$), $\ln(n)^{-2s/\delta}$ ($\delta>0, r=0$) and
$\ln(n)^{(2\gamma+1)/r}/n$ ($ \delta=0, r>0$) are known to be the
minimax rates and we refer to Fan~(1991) (first two cases) and to
Butucea~(2004) (last case) for lower bounds.
\\

\noindent The optimality of the rates in the case $\delta>0, r >0$ requires a
specific discussion.

To our knowledge, the first paper dealing with the case where $g$ is super smooth ($r>0$) is the paper by Pensky and Vidakovic~(1999). See Section 4.3 for a discussion of the rates they obtain
compared to ours.

The case $r=\delta=1$ is studied by Tsybakov~(2000) and  Cavalier
et al.~(2003), in the case of inverse problems with random noise.
In this case and in both problems (density deconvolution and
inverse problem) the best compromise is explicit and so is the
rate of convergence, of order
$n^{-a'/(a'+\mu\sigma)}[\ln{n}]^{(-2s
\mu\sigma+2b\gamma)/(\mu\sigma+b)}$. It is noteworthy that $\hat
g_m^{(n)}$ seems also to achieve the minimax rate of convergence
in this case.\newline

When $0<r <\delta$, some lower bounds are known in the special
case $0<r <\delta$ and $s=0$. According to Butucea and
Tsybakov~(2004), in this case, if we denote by $\pi L_{\breve m}$
the solution of $2\mu\sigma^\delta(\pi L_{\breve m})^\delta+2b
(\pi L_{\breve m})^r=\ln{n}-(\ln\ln{n})^2,$ then the rate of
convergence of $\hat g_m$ is the minimax rate of order $\exp\{-2b
(\pi L_{\breve m})^r\}$. The rate of convergence is always of
order a power of $\ln(n)$ multiplied by an exponential term, that
is decreases faster that any logarithmic function, but slower than
any power of $n$.\newline

When $0<\delta< r$, no lower bounds are available. In this case,
the rate is of order a power of $\ln(n)$ multiplied by a negative
power of $n$ and by an exponential term.

\bigskip\noindent\textit{3.5 Conclusion on the minimum contrast estimators
$\hat g_m^{(n)}$}

The estimator $\hat g_m^{(n)}$ achieves the
minimax rate in all cases where lower bounds are available but its
construction requires the knowledge of the smoothness of $g$. All
those facts give strong motivation to  find some adaptive
estimation procedure that does not require such prior smoothness
knowledge on $g$, and whose risk automatically achieves the
minimax rate.

\bigskip\noindent\textsf{4.\ ADAPTIVE ESTIMATION}

\bigskip\noindent\textit{4.1 Main result of adaptive estimation}

\bigskip\noindent We look for a penalty function, based on the observations
and on $\sigma f_\varepsilon(\cdot/\sigma)$, such that, for $K_n\geq n$
\begin{equation}
\label{oracleind} {\mathbb E}\parallel \tilde g -g\parallel^2 \leq
\inf_{m \in \mathcal{M}_n}\left[\parallel g-g_m\parallel^2+(\pi L_m)^2(M_2+1)/n+
2\lambda_1 \Gamma(m)/n\right].
\end{equation}
The following theorem describes the cases where the oracle
inequality (\ref{oracleind}) is reached.\newline

\noindent\textsc{Theorem 1}. \textit{
Under the assumptions {\rm (A}$_2^{\varepsilon}${\rm )} and {\rm (A}$_3^X${\rm )}, consider the
collection of estimators $\hat g_m^{(n)}$ defined by
(\ref{tronque})
with $K_n\geq n$ and $1\leq m\leq m_n$  satisfying (\ref{mn}) if $\delta\leq 1/3$ and if $\delta>1/3$,
\begin{eqnarray*}
m_n\leq
\pi^{-1}\left[\frac{\ln(n)}{2\mu\sigma^\delta}+\frac{2\gamma+1-\delta+
\min((3\delta/2-1/2),\delta)}{2\delta\mu
\sigma^\delta}\ln\left(\frac{\ln(n)}{2\mu\sigma^\delta}\right)\right]^{1/
\delta}.
\end{eqnarray*}
Let {\rm pen}$(m)$ be defined by (\ref{penalite}) for some universal numerical constant $a>1$. Then, $\tilde g = \hat g_{\hat m}^{(n)}$ defined by (\ref{estitronc})
satisfies
\begin{equation}\label{resu1}
{\mathbb E}(\|g-\tilde g\|^2) \leq C_a\inf_{m\in \{1, \dots,
m_n\}}[\|g-g_m\|^2+ {\rm
pen}(m)+(\pi L_m)^2(M_2+1)/n] +
a\kappa_a C/n,
\end{equation} where $C_a=\max(\kappa_a^2, 2\kappa_a)$, $\kappa_a=(a+1)/(a-1)$
and $C$ is a constant depending on $f_{\varepsilon}$ and $\sigma$.} \newline

Obviously, Remark 1 still holds for the adaptive estimator.

The rates are easy to deduce from (\ref{resu1}) as soon as $g$
belongs to some smoothness class, but the procedure will reach the
rate without requiring the knowledge of any smoothness
parameter.\newline

\bigskip\noindent\textit{4.2 About the optimality of the adaptive estimator
$\tilde g$}

\bigskip \noindent {\textit{ Rate of $\tilde g$ under {\rm (R}$_2^X${\rm
)}~: no loss. }}  \newline

If $g$ satisfies {\rm (R}$_2^X${\rm )},
then according to Section~3.2, $\|g-g_m\|^2=0$ as soon as $\pi
L_m\geq d$, and the parametric rate of convergence is
automatically achieved without the knowledge of $C_2$ and $d$ and
especially without requiring to know that {\rm (R}$_2^X${\rm )} is
fulfilled.

\bigskip\noindent {\textit{ Rate of $\tilde g$ under {\rm (R}$_1^X${\rm
)}.}}\newline

Under {\rm (R}$_1^X${\rm )}, the rate of convergence of $\tilde g$ clearly depends on the order of the
penalty compared to the variance order $\Gamma(m)/n$.  If $g$ satisfies {\rm (R}$_1^X${\rm )}, $\|g-g_m\|^2\leq
(C_1/2\pi)L_m^{-2s}\exp\{-2b \pi^r L_m^r\}$. For instance, if $\delta=0$, by
associating the order of the bias to
the value of pen$(m)$, of order $\Gamma(m)/n$,
we obtain that the estimator $\tilde g$ automatically reaches the minimax rate $\ln(n)^{(2\gamma+1)/r}/n$,
without the knowledge of $s, r$ nor $b$. In all cases, $\tilde g$ achieves the
minimax rate up to some logarithmic factor.

\bigskip\noindent {\textit{ Rate of $\tilde g$ under {\rm (R}$_1^X${\rm )}, cases without loss.}}

\indent When $0\leq \delta\leq 1/3$, the penalty function has the variance
order $\Gamma(m)/n$, and $\tilde g$ achieves the best rate of
$\hat g_{\breve m}$. Under {\rm (R}$_1^X${\rm )}, this best rate
is the minimax rate in all cases here, except if $r\geq \delta>0$
and $\delta\leq 1/3$ which is a case where no lower bounds are
available.

When $\delta>1/3$, the penalty function pen$(m)$ has not exactly
the order of the variance $\Gamma(m)/n$, but a loss of order
$L_m^{\min((3\delta/2-1/2), \delta)}$ occurs, that is of order
$L_m^{(3\delta-1)/2}$ if $1/3<\delta\leq 1$ and of order
$L_m^{\delta}$ if $\delta>1$. Consequently $\tilde g$ achieves the
best rate of $\hat g_{\breve m}$ if the bias $\|g-g_m\|^2$ is the
dominating term in the trade-off between $\|g-g_m\|^2$ and
$\mbox{pen}(m)$.

\indent \textbf{-} When $r=0$ and $\delta>1/3$, the minimax rate of
order $(\ln(n))^{-2s/\delta}$ is given by the bias term, and
the loss in the penalty function does not change the rate achieved
by the adaptive estimator $\tilde g$, which remains thus the minimax rate.

\indent \textbf{-} When $0<r<\delta$, the rate
is given by the bias term and thus this loss does
not affect the rate of convergence of $\tilde g$ either.
Therefore, $\tilde g$ achieves the best rate of $\hat g_{\breve m}$, which is
the minimax rate of convergence when $s=0$ and also probably if $s\not=0$. In the
specific case $0<r<\delta/2$ and $s=0$, Butucea and Tsybakov~(2004) also propose an
adaptive estimator. But this requires to know that $0<r<\delta/2$ and $s=0$.

\bigskip\noindent {\textit{ Rate of $\tilde g$ under {\rm (R}$_1^X${\rm )}, case with loss.}}\newline
\indent \textbf{-} When $r\geq \delta>1/3$,  $\mbox{pen}(m)$ can
be the dominating term in the trade-off between $\|g-g_m\|^2$ and
$\mbox{pen}(m)$. This induces a loss of order
$L_m^{\min((3\delta/2-1/2),\delta)}$ in the rate of convergence of
$\tilde g$ compared to the best rate of $\hat g_{\breve m}$. Since
it happens in cases where the order of the optimal $L_m$ is less
than $(\ln{n})^{1/\delta}$, the loss in the rate is at most of
order $\ln{n}$, when the rate is faster than logarithmic and
consequently, the loss appears only in cases where it can be seen
as negligible.

For $\mathbb{L}_2$ estimation, such an unavoidable logarithmic
loss in adaptation, has been pointed out by Tsybakov~(2000) and
Cavalier et al.~(2003) in case of inverse problems with random
noise, when $r=\delta=1$, which shows, in a slightly different
model but with comparable rates of convergence, that a loss due to
adaptivity of order $\ln(n)^{b/(\mu\sigma+b)}$ is unavoidable. The
main point is that, according to (\ref{nudelopt}), our estimator
has its quadratic risk with the same logarithmic loss when
$r=\delta=1$. This logarithmic loss due to adaptation seems thus
unavoidable at least in one case.\newline

\noindent\textsc{Remark 2.}
When $\sigma=0$, then by convention $\delta=\mu=0$,
$\lambda_1=1$ and ${\rm pen}(m)=6 a L_m/n$ which is  the penalty
function used in direct density estimation. More precisely, if
$\sigma$ is very small, then the procedure selects the parameter
$L_m$ closed to the parameter selected in usual density
estimation.

\bigskip\noindent\textit{4.3 Comparison with Pensky and
Vidakovic~(1999)}

\bigskip\noindent
To our knowledge, the first paper dealing with adaptive density deconvolution
is the paper by  Pensky and
Vidakovic~(1999) who are also the first that consider the case of $r>0$.
The adaptive estimators proposed in  Pensky and
Vidakovic~(1999) achieve minimax rates of convergence in the three
cases $(\delta=0, r=0)$, $(\delta=0, r>0)$, and $(\delta>0,
r=0)$. 

But when $(r>0, \delta>0)$, the rate of convergence of their
estimator is not minimax. This is shown in the special case
$0<r<\delta$ and $s=0$, in Butucea and Tsybakov~(2004), where
sharp minimax results are stated. This is also shown by our
results when $0<\delta\leq r$ and when $0<r<\delta$, $s\not=0$ (see Sections 3.4 and 4.2). 
For instance, when $0<\delta/r\leq 1/2$, according to (\ref{comparevit}) and Sections 3.3 and 4.2, the resulting rate of $\tilde g$ is
of order
\begin{eqnarray*} 
\ln(n)^{\max ( 0, \min(3\delta/2-1/2,\delta)/r )}\ln(n)^{(2\gamma+1-\delta)/r} \exp[2\mu
\sigma^{\delta}({\ln(n)}/(2b))^{\delta/r}]/n
\end{eqnarray*}
stricly faster than the upper bound of the rate in Pensky
and Vidakovic~(1999) (see their Theorem 4) which is of order
$\ln(n)^{(2\gamma+1)/\delta}/n^{\Lambda/(\Lambda+2\mu\sigma^\delta
(4\pi/3)^\delta)}$ for $\Lambda>0$. 

The non-optimality of their adaptive estimator when $(\delta>0, r>0)$ comes from
two facts. First, when $(\delta>0, r>0)$, they choose a smoothing
parameter (analogous to $L_{\breve m}$) as in the case
$(r=0,\delta>0)$. Consequently, it provides an adaptive estimator in the
sense that it does not depend on the smoothness parameters of $g$. But it 
does not
give the best rate for their estimator, since it does not correspond to the best choice in
their bias-variance compromise. 

Second, this non optimality of their estimator when $\delta>0,~r>0$, comes also, in a more crucial manner, from the fact that
their wavelet and scaling functions cannot provide the optimal bias-variance
decomposition. This is due to the
support of the Fourier transform of their scaling function  as
well as their wavelet which induce, when $\delta>0,~r>0$, a squared bias term of order $L_m^{-2s}\exp\{-2b
(2\pi/3)^{r}L_m^r \}$ with a variance term of order
$L_m^{2\gamma+1-\delta}\exp\{2\mu\sigma^\delta
(4\pi/3)^{\delta}L_m^\delta\}.$ 
When either $(\delta=0, r=0)$, $(\delta>0, r>0)$ or
$(\delta>0,r=0)$, those supports have no influence on the rate of convergence,
and hence their estimator
is minimax. But these supports do not allow to reach the minimax rate when
$(\delta>0, r>0)$.

The asymptotic properties of $\tilde g$ are 
improved by using the basis generated by 
$\sin(\pi
x)/(\pi x)$. Indeed, due to its Fourier transform, it
implies a squared bias of order
$L_m^{-2s}\exp\{-2b \pi^{r}L_m^r\}$ and a variance of order
$L_m^{2\gamma+1-\delta}\exp\{2\mu\sigma^\delta\pi^\delta
L_m^\delta \}$ and hence a better trade-off between the two terms. Section 3.3 as well as Butucea and
Tsybakov~(2004)'s results illustrate that the best choice of $L_{\breve m}$,
solution of the bias-variance compromise (see equation (\ref{implicit})), requires quite precise computations. 
Besides its simplicity, this basis seems thus the most
relevant since it gives the minimax rates in
all the cases where lower bounds are available and faster rates than the ones
in Pensky and Vidakovic~(1999) in the
remainder case. 


\bigskip\noindent\textsf{5.\ SIMULATION STUDY}

\bigskip\noindent The implementation is conducted by using Matlab
software. Details about the algorithm can obtained from the authors upon
request. We choose $K_n=2^8$ as being of order
$O(n)$ is all cases.

The integrated squared error ISE$(\hat g_{\hat m}^{(n)})= \|\hat g_{\hat m}^{(n)}-g\|^2$ is computed
via a standard approximation and discretization of
the integral on an interval of ${\mathbb R}$ denoted by $I$ and given in each case.

Then the MISE, ${\rm MISE}(\hat g_{\hat m}^{(n)})=\mathbb{E}\|\hat g_{\hat m}^{(n)}-g\|^2$ is
computed as the empirical mean of the approximated ISE $\|\hat
g_{m}^{(n)}-g\|^2$, over 500 simulation samples.
We illustrate our method on some test densities, with various
smoothness properties, and for the two types of errors, ordinary
and super smooth. We start by describing the error densities and
the associated penalties.

\bigskip\noindent\textit{5.1 Two settings for the errors and the associated penalties}

\bigskip\noindent We consider two types of error density $f_\varepsilon$, the first
one is ordinary smooth, with polynomial decay of the Fourier
Transform, and the second one is supersmooth, with an exponential
decay of the Fourier transform $f_\varepsilon^*$.\newline

$\bullet$ \textbf{Case 1: Laplace (or Double exponential)
$\varepsilon$'s.} In this case,
$f_{\varepsilon}(x)= e^{-\sqrt{2}|x|}/{\sqrt{2}},$ and $
f^*_{\varepsilon}(x)= (1+ x^2/2)^{-1}.$

This density satisfies {\rm (A}$_2^{\varepsilon}${\rm )} with $\gamma=2$, $\kappa_0=1/2$ and
$\mu=\delta=0$.\newline
According to Theorem~1, the penalty function, as the
variance, is of order
$$\frac{L_m}{n}\int_{-\pi}^{\pi} \left|\frac{\varphi^*(x)}{f_{\varepsilon}^*(\sigma
{L_m}x)}\right|^2dx, \mbox{ where }
\int_{-\pi}^{\pi} \left|\frac{\varphi^*(x)}{f_{\varepsilon}^*(\sigma {L_m}x)}\right|^2dx =
2\pi  \left(1+ \frac{\pi^2}3 \sigma^2{L_m}^2 + \frac{\pi^4}{20}
\sigma^4{L_m}^4\right).$$
Some intensive simulation studies on various tested densities lead to choose the following penalty
$${\rm pen}({L_m}) = \frac{6\pi {L_m}}n\left( 1+ \frac{(\ln({L_m}))^{2.5}}{{L_m}} + \frac{\pi^2}3
\sigma^2{L_m}^2 + \frac{\pi^4}{20} \sigma^4{L_m}^4\right).$$

$\bullet$ \textbf{Case 2: Gaussian $\varepsilon$'s.} In that case,
$f_{\varepsilon}(x)=1/\sqrt{2\pi}e^{-x^2/2},$ and $
f_{\varepsilon}^*(x)= e^{-x^2/2}.$

This density
satisfies {\rm (A}$_2^{\varepsilon}${\rm )} with $\gamma=0$, $\kappa_0=1$,
$\delta=2$
and $\mu=1/2$.\newline
According to Theorem~1, the penalty, slightly bigger than the variance
term, is of order
$$\frac{{L_m}^{3}}n \int_{-\pi}^{\pi} \left|\frac{\varphi^*(x)}{f_{\varepsilon}^*(\sigma
{L_m}x)}\right|^2dx \quad \mbox{ where }\quad \int_{-\pi}^{\pi}
\left|\frac{\varphi^*(x)}{f_{\varepsilon}(\sigma {L_m}x)}\right|^2dx =
\int_{-\pi}^{\pi} \exp(\sigma^2 {L_m}^2 x^2) dx.$$ As in the
previous case, some intensive simulation studies on various tested densities lead to choose the following penalty
$${\rm pen}({L_m})= \frac{6\pi {L_m}}n \left(1+\frac{(\ln({L_m}))^{2.5}}{L_m} +
  \frac{\pi^2\sigma^2{L_m}^2}3\right)\left(\int_{0}^{\pi} \exp(\sigma^2 {L_m}^2 x^2) dx/\pi\right),$$
where the integral is numerically computed. According to the
theory (see Theorem~1, the loss due to the
adaptation is the term $\pi^2\sigma^2{L_m}^2/3$.\newline
The
additional term $(\ln({L_m}))^{2.5}/{L_m}$ is motivated by the works
of Birg\'e and Rozenholc~(2005). In our case also, this term improves the quality of the
results by making the penalties slightly heavier when ${L_m}$ becomes large. \newline
Note that when $\sigma=0$, both penalties are equal to $(6\pi
{L_m})(1+(\ln({L_m}))^{2.5}/{L_m})/n.$


\bigskip\noindent\textit{5.2 Test densities}

\bigskip\noindent First we consider densities having classical smoothness properties like H\"olderian
smoothness with polynomial decay of their Fourier transform. Second we consider
densities having stronger smoothness properties, with exponential decay of the
Fourier transform.
Except in the case of the infinite variance density (Cauchy
density), we consider density functions $g$ normalized with unit
variance so that $1/\sigma^2$ represents the usual signal-to-noise
ratio (variance of the signal divided by the variance of the noise)
and is denoted in the sequel by $s2n$ defined as $s2n =
1/{\sigma^2}.$
The functions
which are considered are listed below, associated with the interval
$I$ used to evaluate the ISE:\newline
(a) Chi2(3)-type distribution, $X=1/\sqrt{6}U$, $g_X(x)=\sqrt{6}g(\sqrt{6}x)$,
$U\sim \chi^2(3)$ where we know that $U\sim \Gamma(\frac 32, \frac 12)$,

and $I=[-1,16]$.\newline
(b) Laplace distribution, $I=[-5,5]$.\newline
(c) Mixed Gamma distribution, $X=1/\sqrt{5.48} W$ with $W\sim 0.4
\Gamma(5,1) + 0.6 \Gamma(13,1)$,

and $I=[-1.5,26]$.\newline
(d) Cauchy distribution, $g(x)=(1/\pi)(1/(1+x^2))$, $g^*(x)=e^{-|x|}$, $I=[-10,10]$.\newline
(e) Gaussian distribution, $X\sim {\mathcal N}(0,\sigma^2)$ with $\sigma=1$, $I=[-4,4]$.\newline
(f) Mixed Gaussian distribution: $X\sim \sqrt{2}V$ with $V\sim 0.5{\mathcal N}(-3,1)
+ 0.5 {\mathcal N}(2,1)$
and $I=[-8,7]$.\newline

Densities (a), (b), (c) correspond to cases with $r=0$, whereas
densities (d), (e), (f) correspond to cases with $r>0$.

\bigskip\noindent\textit{5.3 Results}

\begin{figure}[htb]
\hspace{2cm} $$\includegraphics[width=0.4\textwidth,
height=15cm,angle=-90 ]{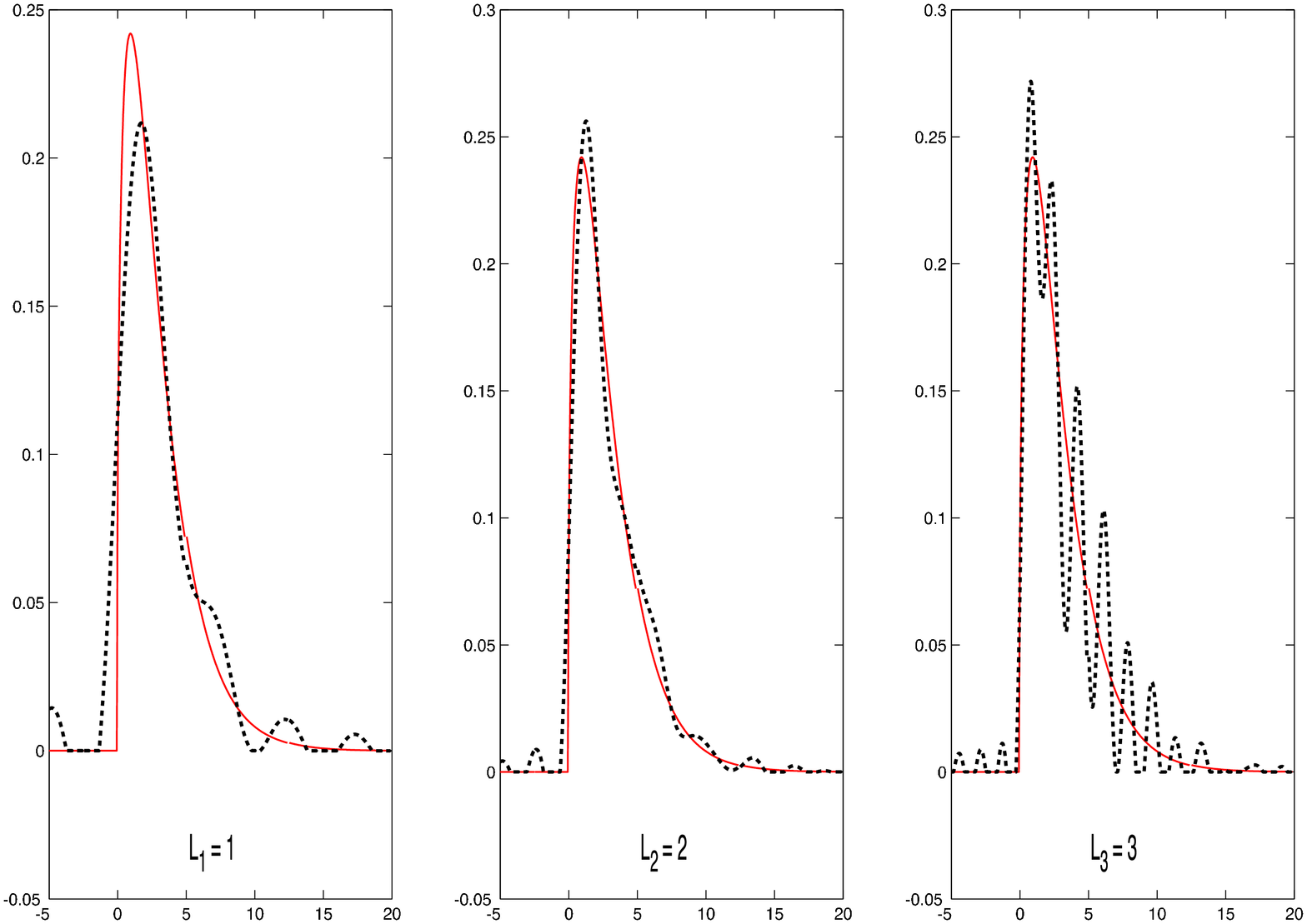}$$
\caption{Plots of the estimator (dotted line) and of the true $\chi^2(3)$ (a)
density  (full line) - Laplace errors -
$n=750$, $s2n$=10, when $L_1=1$ (left), $L_2=2$ (middle), $L_3=3$ (right). The algorithm
chooses $\hat m=L_{\hat m}=2$.}\label{ExChi2Laplace}
\end{figure}

Figure \ref{ExChi2Laplace} compares the estimators $\hat g^{(n)}_m$ obtained
for $m=L_m=1, 2$ and 3, and justifies the good choice $\hat m=2$ of the algorithm.
Table \ref{BasicMise} presents the MISE for the two types of
errors, the different tested densities, different $s2n$ and
for different sample sizes. The greatest values of $s2n$ amount to
consider that there is essentially no noise. Clearly the MISE are
smaller when there is less noise ($\sigma$ small, $s2n$ large).

\begin{table}[ptbh]
\caption{Mean MISE $\times 100$ obtained with $N=500$ samples, for
sample  sizes $n=100, 250, 500, 1000, 2500$ and $s2n=$ 2, 4, 10,
100, 1000, the higher $s2n$ the lower the noise level. Densities
(a): Chi2(3), (b): Laplace, (c): Mixed Gamma, (d): Cauchy, (e)
Gaussian, (f): Mixed Gaussian.}\label{BasicMise}
\begin{center}
{\small
\begin{tabular}{crcccccccccc}\hline
\multicolumn{2}{c}{$\times10^{-2}$}  &
\multicolumn{2}{c}{$n=100$}& \multicolumn{2}{c}{$n=250$}&
\multicolumn{2}{c}{$n=500$}& \multicolumn{2}{c}{$n=1000$}&
\multicolumn{2}{c}{$n=2500$}\\ \hline $g$ & $s2n$ & Lap. &
Gaus. & Lap. & Gaus. & Lap. & Gaus. & Lap. & Gaus. & Lap. & Gaus. \\
\hline  \multirow{5}{.5cm}{(a)}
 &    2  & 2.02 & 4.15 &  1.39 &  2.37 &  1.18 &  1.72 &  1.06 &  1.36 &  1.03 &  1.12\\
 &    4 & 1.52 & 1.79 &  1.21 &  1.27 &  1.07 &  1.13 &  1.04 &  1.04 & 0.654 & 0.996\\
 &   10 & 1.31 & 1.31 &  1.13 &  1.11 &  1.01 &  1.03 & 0.505 & 0.995 & 0.345 & 0.974\\
 &  $10^2$ & 1.22 & 1.23 &  0.72 & 0.884 & 0.409 & 0.411 & 0.327 & 0.335 & 0.179 & 0.232\\
 & $10^3$  & 1.22 & 1.21 & 0.651 & 0.638 & 0.391 & 0.382 & 0.293 & 0.298 & 0.157 & 0.157\\
\hline \multirow{5}{.5cm}{(b)}
 &    2 &  3.7 & 10.6 &  2.17 &  5.2 &  1.61 &  3.03 &  1.41 &  2.07 &   1.2 &  1.48\\
 &    4 &  2.5 & 2.99 &  1.66 & 1.93 &  1.33 &  1.46 &  1.26 &  1.25 & 0.817 &  1.12\\
 &   10 &  1.9 & 1.97 &  1.43 & 1.42 &  1.35 &  1.22 & 0.723 &  1.12 & 0.441 &  1.06\\
 &  $10^2$  & 1.69 & 1.64 & 0.883 & 1.06 & 0.607 & 0.538 & 0.453 & 0.385 & 0.343 & 0.211\\
 & $10^3$   & 1.68 & 1.65 & 0.814 & 0.79 & 0.593 & 0.561 & 0.411 & 0.379 & 0.284 &  0.24\\
\hline \multirow{5}{.5cm}{(c)}
 &    2 &  1.32 &  3.96 & 0.547 &  1.88 &  0.292 &   1.01 &  0.148 &  0.533 &   0.06 &  0.224\\
 &    4 &  0.79 &  1.05 & 0.316 & 0.453 &  0.151 &  0.224 & 0.0815 &  0.116 & 0.0361 & 0.0497\\
 &   10 &  0.495 & 0.524 & 0.194 & 0.215 &  0.103 &   0.11 & 0.0543 & 0.0565 &  0.024 & 0.0246\\
 &  $10^2$ & 0.369 & 0.384 & 0.152 & 0.149 & 0.0789 & 0.0785 & 0.0409 & 0.0412 & 0.0194 & 0.0186\\
 & $10^3$  & 0.364 & 0.353 & 0.149 &  0.15 & 0.0762 & 0.0767 & 0.0404 & 0.0406 & 0.0184 & 0.0185\\
\hline \multirow{5}{.5cm}{(d)}
 &    2  &  2.72 &  9.09 &  1.22 &  4.26 & 0.645 &   2.3 & 0.353 &  1.25 & 0.158 &  0.513\\
 &    4  &  1.66 &  2.27 & 0.716 & 0.967 & 0.364 & 0.514 & 0.205 &  0.28 & 0.138 &  0.127\\
 &   10  &  1.15 &  1.13 & 0.437 &  0.46 & 0.249 & 0.257 & 0.215 & 0.142 & 0.219 & 0.0764\\
 &  $10^2$ & 0.815 & 0.783 & 0.373 & 0.351 & 0.351 & 0.271 & 0.206 & 0.201 & 0.147 & 0.0962\\
 & $10^3$  & 0.783 &  0.78 & 0.366 & 0.355 &  0.34 & 0.331 & 0.189 & 0.189 & 0.121 &  0.118\\
\hline \multirow{5}{.5cm}{(e)}
 &    2  &  2.74 &  9.21 &   1.1 &  4.08 & 0.605 &  2.14 & 0.296 &  1.06 & 0.143 &  0.446\\
 &    4  &  1.59 &  2.23 & 0.591 & 0.878 & 0.362 & 0.457 & 0.229 & 0.227 & 0.463 & 0.0894\\
 &   10  & 0.885 &  1.02 & 0.397 &  0.42 & 0.372 &  0.21 & 0.515 & 0.112 & 0.229 &  0.046\\
 &  $10^2$  & 0.711 & 0.713 & 0.565 & 0.432 & 0.396 & 0.394 & 0.279 & 0.195 & 0.171 &   0.15\\
 & $10^3$   & 0.739 & 0.705 & 0.606 & 0.592 & 0.352 & 0.355 & 0.259 & 0.246 & 0.167 &  0.145\\
\hline \multirow{5}{.5cm}{(f)}
 &    2 &  2.97 &  9.98 &  1.26 &  4.45 & 0.693 &  2.31 &  0.328 &   1.26 &  0.132 &  0.509\\
 &    4 &  1.73 &  2.37 & 0.709 &  1.02 & 0.375 & 0.478 &  0.185 &  0.257 & 0.0751 &  0.105\\
 &   10 &  1.14 &  1.21 & 0.463 & 0.466 & 0.237 & 0.242 &  0.118 &  0.122 & 0.0468 & 0.0515\\
 &  $10^2$ & 0.851 & 0.817 & 0.359 & 0.352 & 0.166 & 0.167 & 0.0866 & 0.0867 &  0.034 & 0.0351\\
 & $10^3$  & 0.823 & 0.828 & 0.344 & 0.327 & 0.169 & 0.163 & 0.0845 & 0.0839 & 0.0334 & 0.0336\\
 \hline
\end{tabular}}
\end{center}
\end{table}

We can in particular compare the performances of our adaptive
estimator with the performances of the deconvolution kernel as
presented in Delaigle and Gijbels~(2004a). This comparison  is done
for densities (a), (c), (e) and (f) which correspond to the
densities $\#2$, $\#6$, $\#1$ and $\#3$ respectively, in  Delaigle
and Gijbels~(2004a). They give median ISE obtained with kernel
estimators by using four different methods of bandwidth selection.
The comparison is given in Table \ref{DelGibcomp} between the median
ISE computed for 500 samples generated with the same interval length
and signal to noise ratio as Delaigle and Gijbels (2004a). The ISE
are computed on the same intervals $I$ as them. We  also give our
corresponding means since we believe that they are more meaningful
than medians since the MISE is $\mathbb{E}\|\hat g_{m}^{(n)}-g\|^2$,
but we also give our medians.

We can see that our estimation procedure provides better results in
all cases except in one case, namely when we aim at estimating a
Gaussian density, for both types of errors density. This is the most
probably due to the fact that the bandwidth selection methods are
based on computations assuming that the underlying density is
Gaussian, so that they perform very well when it is true. For the
other cases, even our means are often better than Delaigle and
Gijbels'(2004a) medians which shows that our method provides a very
good solution to the deconvolution problem.

A standard objection to deconvolution methods is that they require
the knowledge of the noise density. Therefore, following the ideas
of Meister~(2004), we study here the properties of the estimator
when the error density is not correctly specified. For both type
of errors, we study the behavior of the estimator using one type
of the error density when the other type of errors density is the
good one. Table \ref{errorinerror} presents the ratio between the
resulting MISE if the errors density is not correct with the MISE
if the errors density is correct. For instance, in the columns
"$\varepsilon$ Lap." the noise density is Laplace but the MISE in
the numerator of the ratio corresponds to estimators constructed
as if it were Gaussian. As expected, since the construction uses
the knowledge of the error density, if it is misspecified, the
estimator presents some bias and the MISE becomes slightly bigger.
Nevertheless, this difference does not clearly appear when $n$ is
not very large. Indeed in that case, the optimal length ${L_m}$ is
small and therefore the variance term of order $\int_0^{\pi {L_m}}
\vert f_\varepsilon^*(x)\vert^{-2}dx$ is not so different between
the two errors.

\begin{table}[ptbh]
\caption{Median ISE obtained by Delaigle and Gijbels (2004a) with a
kernel estimator and four different strategies of bandwidth
selection, and with our penalized projection estimator (median and
mean).}\label{DelGibcomp}
\begin{center}
{\small
\begin{tabular}{clcccc}\hline
\multicolumn{2}{c}{$\;$} & \multicolumn{2}{c}{$n=100$} &
\multicolumn{2}{c}{$n=250$}\\\hline  density $g$ & method &
$\varepsilon$ Lap. & $\varepsilon$ Gaus. & $\varepsilon$ Lap. &
$\varepsilon$ Gaus. \\ \hline
\multirow{4}{3cm}{\hspace{0.7cm}(a) or $\# 2$\\
\hspace{1cm}$\chi^2(3)$ \\\hspace{0.7cm}($s2n$=4)}
 & DG, lower median& 0.015 & 0.018 & --- & --- \\
 & DG,  higher median  & 0.018 & 0.022 & --- & --- \\ \cline{2-6}
& Proj.: median & 0.014 & 0.016 & --- & --- \\
& Proj.: mean & 0.015 & 0.018 & --- & --- \\  \hline
\multirow{4}{3cm}{\hspace{0.7cm}(c) or $\# 6$\\
\hspace{0.4cm}Mix.Gamma \\\hspace{0.7cm}($s2n$=10)}
&  DG, lower median & --- & --- & 0.0021 & 0.0023 \\
& DG,  higher median  & --- & --- & 0.0024 & 0.0026 \\ \cline{2-6}
& Proj.: median  & --- & --- & 0.0017 & 0.0020 \\
& Proj., mean & --- & --- & 0.0019 & 0.0021\\\hline
\multirow{4}{3cm}{\hspace{0.7cm}(e) or $\# 1$\\ \hspace{0.7cm}
Gauss\\\hspace{0.7cm}($s2n$=4)}
& DG, lower median & 0.0071 & 0.0080 & 0.0041 & 0.0051 \\
& DG, higher median & 0.011 & 0.012 & 0.0059 & 0.0072 \\ \cline{2-6}
& Proj.: median & 0.012 & 0.017 & 0.0049 & 0.0066 \\
& Proj.: mean & 0.016 & 0.022 & 0.0059 & 0.0088\\  \hline
\multirow{4}{3cm}{\hspace{0.7cm}(f) or $\# 3$\\
\hspace{0.6cm}Mix.Gauss\\\hspace{0.7cm}($s2n$=4)}
& DG, lower median &  0.018 & 0.027 & 0.011 & 0.020 \\
& DG,  higher median  & 0.031 & 0.034 & 0.023 & 0.028 \\ \cline{2-6}
& Proj.: median   & 0.016 & 0.022 & 0.0063 & 0.0088 \\
& Proj.: mean & 0.017 & 0.024 & 0.0071 & 0.010\\ \hline
\end{tabular}}
\end{center}
\end{table}

\begin{table}[ptbh]
\caption{Ratio between MISE with misspecified error
density (Laplace errors, $g$ estimated as if errors were Gaussian and reciprocally)
and MISE with correctly specified error density.}\label{errorinerror}
\begin{center}
{\small
\begin{tabular}{crcccccccc}\hline
\multicolumn{2}{c}{$\times 10^{-2}$}  & \multicolumn{2}{c}{$n=1000$}& \multicolumn{2}{c}{$n=5000$}&
\multicolumn{2}{c}{$n=10000$}& \multicolumn{2}{c}{$n=25000$}
\\ \hline
$g$ & $s2n$ & $\varepsilon$ Lap. & $\varepsilon$ Gaus. & $\varepsilon$ Lap. &  $\varepsilon$ Gaus. & $\varepsilon$ Lap. & $\varepsilon$ Gaus. & $\varepsilon$ Lap. & $\varepsilon$ Gaus.   \\ \hline
Lapl.
  & 2 & 1.6 & 1.4 & 2.2 & 1.8 & 2.3 & 2.9 & 2.4 & 4.5\\
  & 4 & 1 & 1.3 & 1 & 1.9 & 1 & 2.2 & 1 & 2.3\\
\hline
Mix.Gam.
  & 2 & 1 & 1.1 & 1.3 & 1.6 & 1.6 & 2.1 & 2.2 & 3\\
  & 4 & 1 & 1   & 1.1 & 1.2 & 1 & 1.3 & 1.1 & 1.5\\
\hline
Cauchy
  & 2 & 1.3 & 1.3 & 1.7 & 1.6 & 2.5 & 1.2 & 3.7 & 1.5\\
  & 4 & 1.1 & 1 & 1.2 & 1.1 & 1.3 & 1.1 & 1.4 & 1.2\\
\hline
Gauss
& 2 &  1.1 & 1.4 & 1.4 & 1.1 & 2 & 1 & 3.1 & 1.2\\
  & 4 & 1 & 0.81 & 1.2 & 1   & 1.2 & 1 & 1.8 & 1.3\\
\hline
\end{tabular}}
\end{center}
\end{table}

\bigskip \textit{Concluding remarks :} Our estimation procedure
provides an adaptive estimator which achieves the minimax rate of
convergence (up to a possible logarithmic factor) in all the cases
where lower bounds are available, without any prior smoothness
knowledge on the unknown density $g$. In particular it solves
almost in the best way the bias-variance problem when the best
compromise would not be easily computable.  Furthermore,
this estimation procedure induces a fast practical  algorithm with pretty good practical
results.

\bigskip\noindent\textsf{6.\ PROOFS}

\bigskip\noindent\textit{6.1 Proof of Proposition 1.}

\bigskip\noindent According to (\ref{tronque}), for any given $ m$ belonging to
$\mathcal{M}_n$, $\hat g_{m}^{(n)}$ satisfies, $\gamma_n(\hat
g_{m}^{(n)} )-\gamma_n(g_{m}^{(n)})\leq 0.$ Denoting by $\nu_n(t)$
the centered empirical process
\begin{eqnarray}\label{nu}
\nu_n(t)=\frac 1n \sum_{i=1}^n \left[u_t^*(Z_i)-\langle t,g\rangle
\right],\end{eqnarray} we have that
\begin{eqnarray}
\label{difgamma}
\gamma_n(t)-\gamma_n(s)=\|t-g\|^2-\|s-g\|^2-2\nu_n(t-s),
\end{eqnarray} and therefore,
$\|g-\hat g_m^{(n)}\|^2\leq  \|g-g_m^{(n)}\|^2 + 2\nu_n(\hat
g_m^{(n)} - g_m^{(n)}).$ Since $\hat a_{m,j}-a_{m,j}=
\nu_n(\varphi_{m,j})$, we get that
$$\nu_n(\hat g_m^{(n)}-g_m^{(n)})=\sum_{\vert j\vert \leq K_n} (\hat
a_{m,j}-a_{m,j})\nu_n(\varphi_{m,j}) = \sum_{\vert j\vert\leq K_n}
[\nu_n(\varphi_{m,j})]^2,$$
 and consequently
${\mathbb E}\|g-\hat g_m^{(n)}\|^2\leq  \|g-g_m^{(n)}\|^2 +
2\sum_{j\in \mathbb{Z}}{\rm
Var}[\nu_n(\varphi_{m,j})].$
Now, since the $X_i$'s and the $\varepsilon_i$'s
are independent and identically distributed random variables, we get that
${\rm Var}[\nu_n(\varphi_{m,j})]=n^{-2} \sum_{i=1}^n {\rm
Var} \left[u_{\varphi_{m,j}}^*(Z_i)\right] = n^{-1} {\rm Var}
\left[u_{\varphi_{m,j}}^*(Z_1)\right].$

Apply Lemma 2 to get that
$\sum_{j\in \mathbb{Z}} {\rm
Var} [\nu_n(\varphi_{m,j})] \leq \Delta_1(m)/n,$
 where $\Delta_1(m)$ is defined in Proposition~1. It remains to study $\|g-g_m^{(n)}\|^2$. By
applying Pythagoras Theorem, we have $\|g-g_m^{(n)}\|^2= \parallel
g-g_m\parallel^2+\|g_m-g_m^{(n)}\|^2,$ where $
\|g_m-g_m^{(n)}\|^2=\sum_{\vert j\vert > K_n}a_{m,j}^2\leq
(\sup_{j}j a_{m,j})^2\sum_{\vert j\vert > K_n }j^{-2}.$ Now we
write that
\begin{eqnarray*}
j a_{m,j} &=&j\sqrt{L_m}\int \varphi(L_mx-j)g(x)dx\\
&\leq & L_m^{3/2}\int \vert x\vert \vert \varphi(L_mx -j)\vert
g(x)dx+\sqrt{L_m}\int \vert L_mx-j\vert \vert \varphi(L_mx -j)\vert
g(x)dx\\
&\leq & L_m^{3/2}\left( \int \vert \varphi(L_mx -j)\vert^2
dx\right)^{1/2}\left( \int x^2g^2(x)dx\right)^{1/2}+\sqrt{L_m}\sup_{x}\vert x\varphi(x)\vert.
\end{eqnarray*}
This implies finally that $j a_{m,j}\leq
L_m(M_2)^{1/2}+\sqrt{L_m},$ and Proposition~1 follows.\hfill$\Box$

\bigskip\noindent\textit{6.2 Proof of Theorem~1}

\bigskip\noindent By definition, $\tilde g$ satisfies that for all $m\in {\mathcal
M}_n$, $\gamma_n(\tilde g)+ {\rm pen}(\hat m)\leq
\gamma_n(g_m^{(n)}) + {\rm
  pen}(m).$
Therefore, by applying (\ref{difgamma}) we get that
\begin{eqnarray*}
\parallel \tilde g-g\parallel^2&\leq& \parallel
g_m^{(n)}-g\parallel^2+2\nu_n(\tilde g -g_m^{(n)})+{\rm pen}(m)-{\rm pen}(\hat m).
\end{eqnarray*}
Next, we use  that if $t=t_1+t_2$ with $t_1$ in $S_m^{(n)}$ and
$t_2$ in $S_{m'}^{(n)}$, then $t$ is such that $t^*$ has its
support in $[-\pi L_{\max(m,m')}, \pi L_{\max(m,m')}]$ and
therefore $t$ belongs to $S_{\max(m,m')}^{(n)}$. If we denote by
$B_{m, m'}(0,1)$ the set $B_{m, m'}(0,1)=\{t\in
S_{\max(m,m')}^{(n)} \;/\; \|t\|=1\},$ then $|\nu_n(\tilde
g-g_m^{(n)}) |\leq \|\tilde g-g_m^{(n)}\|\sup_{t\in
  B_{m,\hat m}(0,1)}|\nu_n(t)|.$
Consequently, by using that $2uv \leq
a^{-1}u^2+av^2$, for $a>1$, we get
\begin{eqnarray*} \|\tilde g-g\|^2 &\leq& \|g_m^{(n)}
-g\|^2 + a^{-1}\|\tilde g-g_m^{(n)}\|^2  + a\sup_{t\in B_{m,\hat
m}(0,1)}\nu_n^2(t)+ {\rm pen}(m)- {\rm pen}(\hat m)
\end{eqnarray*}
and therefore, by writing that $\|\tilde g-g_m^{(n)}\|^2\leq (1+y^{-1})\|\tilde
g-g\|^2+ (1+y)\|g-g_m^{(n)}\|^2$, with $y=(a+1)/(a-1)$ for $a>1$, we infer
that
\begin{eqnarray*} \|\tilde g-g\|^2
\leq \left(\frac{a+1}{a-1}\right)^2 \|g-g_m^{(n)}\|^2 +
\frac{a(a+1)}{a-1}\!\!\!\sup_{t\in B_{m,\hat m}(0,1)}\!\!\!\nu_n^2(t)
+\frac{a+1}{a-1} ({\rm pen}(m)- {\rm pen}(\hat m)).
\end{eqnarray*}
Choose some positive function $p(m,m')$ such that $a p(m,m')\leq {\rm pen}(m)
+ {\rm pen}(m')$. Consequently, for $\kappa_a=(a+1)/(a-1)$ we have
$$\|\tilde g-g\|^2
  \leq \kappa_a^2
  \left[\|g-g_m\|^2+\|g_m-g_m^{(n)}\|^2+{\rm pen}(m)\right] +   a \kappa_a W_n(\hat m)
  $$
\begin{equation}\label{Wg} \mbox{ with }\hspace{1.5cm}W_n(m'):=[\sup_{t\in B_{m, m'}(0,1)}
 |\nu_n(t)|^2-p(m, m')]_+,
\end{equation}
that is, according to the proof of Proposition~1,
\begin{equation}\label{majo2} \|\tilde g-g\|^2
  \leq \kappa_a^2 \|g-g_m\|^2+\kappa_a^2(M_2+1)(\pi L_m)^2/K_n
  + 2\kappa_a{\rm pen}(m) + a\kappa_a\sum_{m'\in
    \mathcal{M}_n}W_n(m').\end{equation} The main point of the
proof lies in studying $W_n(m')$,
and more precisely in finding $p(m,m')$ such that for a constant $K$,
\begin{equation}\label{but}
\sum_{m'\in {\mathcal M}_n} \mathbb{E}(W_n(m'))\leq K/n.
\end{equation}
In this case, combining (\ref{majo2}) and (\ref{but}) we infer that, for all
$m$ in ${\mathcal M}_n$,
$$  \mathbb{E}\|g-\tilde g\|^2 \leq \kappa_a^2
  \|g-g_m\|^2 + \kappa_a^2 (M_2+1)(\pi L_m)^2/K_n
  +2\kappa_a {\rm pen}(m)+ a\kappa_a K/n,$$
 which can also be written
$$\mathbb{E}\|g-\tilde g\|^2 \leq C_a\inf_{m\in {\mathcal M}_n}
\left[ \|g-g_m\|^2 + {\rm pen}(m)+(M_2+1)(\pi L_m)^2/K_n\right] +
a\kappa_a K/n,$$
where $C_a=\max(\kappa_a^2, 2\kappa_a)$ suits. It remains thus to
find $p(m,m')$ such that (\ref{but}) holds. This will be done by
applying the following immediate integration of Talagrand's
Inequality (see Talagrand (1996)):\newline

\noindent\textsc{Lemma 1}. {\textit Let $Y_1, \dots, Y_n$ be
i.i.d. random variables and  $r_n(f)=(1/n)\sum_{i=1}^n
[f(Y_i)-\mathbb{E}(f(Y_i))]$ for $f$ belonging to a countable
class ${\mathcal F}$ of uniformly bounded measurable functions.
Then for $\xi^2>0$
\begin{equation}
\label{talesp} \mathbb{E}\left[\sup_{f\in {\mathcal
F}}|r_n(f)|^2-2(1+2\xi^2)H^2\right]_+ \leq \frac 6{K_1}\left(\frac
vn e^{-K_1\xi^2 \frac{nH^2}v} + \frac{8M_1^2}{K_1n^2C^2(\xi^2)}
e^{-\frac{K_1 C(\xi)\xi}{\sqrt{2}}\frac{nH}{M_1}}\right),
\end{equation}
with $C(\xi)=\sqrt{1+\xi^2}-1$, $K_1$ is a
universal constant, and where $$\sup_{f\in {\mathcal
F}}\|f\|_{\infty}\leq M_1, \;\;\;\; \mathbb{E}[\sup_{f\in
{\mathcal F}}|r_n(f)|]\leq H, \;\;\;\; \sup_{f\in {\mathcal
F}}{\rm Var}(f(Y_1)) \leq v.$$} \newline

Usual density arguments show that this result can be applied to the class of
functions ${\mathcal F}= B_{m,m'}(0,1)$.
 Let us denote by $m^*=\max(m,m')$.  Combining Lemma 3 and Lemma 4, we propose to take
$$H^2= H^2(m^*)= \lambda_1
L_{m^*}^{2\gamma+1-\delta}\exp\{2\mu\sigma^\delta(\pi
L_{m^*})^{\delta}\}/n \mbox{ and } \;\; M_1= \sqrt{nH^2},$$ where
$\lambda_1=\lambda_1(\gamma,\kappa_0,\mu,\sigma,\delta)$ is
defined by (\ref{lambda1}).  Again, by applying Lemma 4, we take
$v\geq \Delta_2(m^*,h)$ with
\begin{equation}\label{Delta2}
\Delta_2(m, h)=  L_m^2\iint \left\vert
\frac{\varphi^*(x)\varphi^*(y)}{f_\varepsilon^*(\sigma L_mx)f_\varepsilon^*(\sigma L_my)}h^*(L_m(x-y))
 \right\vert ^2dxdy.
 \end{equation}
 For $\delta>1$ we use a rough bound for $\Delta_2(m,h)$ given by
 $\sqrt{\Delta_2(m^*,h)}\leq 2\pi nH^2$. When $\delta \leq 1$, write that
\begin{eqnarray*} \Delta_2(m,h) \!\!\!&\leq &\!\!\! \kappa_0^{-2}L_m^2
 (1+(\sigma\pi L_m)^2)^{\gamma} \exp\{2\mu\sigma^\delta(\pi L_m)^{\delta}\}
\int_{-\pi}^{\pi}\frac{dx}{|f_{\varepsilon}^*(\sigma L_m x)|^2} \int
|h^*(L_mu)|^2du \\ \!\!\!&\leq & \!\!\!2\kappa_0^{-2}\pi
\lambda_1(1+\sigma^2\pi^2)^{\gamma}\|h^*\|^2L_{m}^{4\gamma+1-\delta}\exp\{4\mu\sigma^\delta(\pi
L_{m})^{\delta}\}.
\end{eqnarray*}
Using that $\|h^*\|^2 \leq \|f_\varepsilon^*\|^2<\infty$, we take $v= \lambda_2 L_{m^*}^{2\gamma+
  \min (1/2-\delta/2,1-\delta)}\exp\{2\mu\sigma^\delta(\pi
L_{m^*})^{\delta}\},$ where
$\lambda_2=\lambda_2(\gamma,\kappa_0,\mu,\sigma,\delta)$ is
defined in Theorem~1.  From the definition
(\ref{Wg}) of $W_n(m')$, by taking $p(m,m')=2(1+2\xi^2)H^2$, we
get that
$$\mathbb{E}(W_n(m'))\leq\mathbb{E}[\sup_{t\in B_{m,m'}(0,1)}
  |\nu_n(t)|^2 -2(1+2\xi^2)H^2]_+.$$
By applying (\ref{talesp}), we get the global bound
$\mathbb{E}(W_n(L_{m'}))\leq K[I(L_{m^*})+II(m^*)],$ where
$I(m^*)$ and $II(m^*)$ are defined by
\begin{eqnarray*}
&&I(m^*)=\frac{\lambda_2L_{m^*}^{2\gamma+
\min(1/2-\delta/2,1-\delta)}\exp\{2\mu\sigma^\delta(\pi
L_{m^*})^{\delta}\}}{n}
\exp\{{-K_1\xi^2(\lambda_1/\lambda_2)L_{m^*}^{(1/2-\delta/2)_+}}\}\\
\mbox{ and } &&II(m^*)= \frac{\lambda_1
L_{m^*}^{2\gamma+1-\delta}e^{2\mu\sigma^\delta(\pi L_{m^*})^{\delta }}}{n^2}
\exp\left\{-K_1\xi C(\xi)\sqrt{n}/\sqrt{2}\right\},
\end{eqnarray*}
with $\lambda_2=\lambda_2(\gamma,\kappa_0,\mu,\sigma,\delta)$
defined in Theorem~1.

\noindent $\bullet$ Study of $\sum_{m'\in \mathcal{M}_n}II(m^*)$.
We have
$\sum_{m'\in \mathcal{M}_n}II(m^*)\leq \vert \mathcal{M}_n\vert \exp\left\{-K_1\xi C(\xi)
  \sqrt{n}/\sqrt{2}\right\}2\lambda_1\Gamma(m_n)/n^2,
$ according to the choices for $v$, $H^2$ and $M_1$.
Consequently, since under (\ref{mn}), $\Gamma(m_n)/n$ is bounded,
$\sum_{m'\in \mathcal{M}_n}II(m^*)\leq C/n$.

\noindent $\bullet$ Study of $\sum_{m'\in \mathcal{M}_n}I(m^*)$.
Denote by
$\psi=2\gamma+ \min(1/2-\delta/2,1-\delta)$, $\omega=(1/2-\delta/2)_+$,
$K'=K_1\lambda_1/\lambda_2$, then for $a,b\geq 1$, we infer that
\begin{eqnarray}\nonumber
\max(a,b)^{\psi}e^{2\mu\sigma^\delta\pi^{\delta} \max(a,b)^{\delta}}e^{-K'\xi^2\max(a,b)^{\omega}}
&\leq&
(a^{\psi}e^{2\mu\sigma^\delta\pi^{\delta} a^{\delta}}+b^{\psi} e^{2\mu\sigma^\delta\pi^{\delta}
b^{\delta}})e^{-(K'\xi^2/2)(a^{\omega} +
b^{\omega})}\\ \hspace{1cm}\leq   a^{\psi}e^{2\mu\sigma^\delta\pi^{\delta}
a^{\delta}}e^{-(K'\xi^2/2)a^{\omega}} e^{-(K'\xi^2/2)
b^{\omega}}&+&b^{\psi} e^{2\mu\sigma^\delta\pi^{\delta} b^{\delta}} e^{-(K'\xi^2/2)
b^{\omega}}.\label{eqmax}\end{eqnarray}
Consequently, if we denote by $\tilde \Gamma$ the quantity
$\tilde \Gamma(m)=L_{m}^{2\gamma+
\min(1/2-\delta/2,1-\delta)}\exp\{2\mu\sigma^\delta(\pi
L_{m})^{\delta}\}$
then
\begin{eqnarray}
\sum_{m' \in \mathcal{M}_n}I(m^*)
&\leq& \frac{2\lambda_2\tilde \Gamma(m)}{n}\exp\{-(K'\xi^2/2)({L_m})^{(1/2-\delta/2)_+}\}\sum_{m' \in
  \mathcal{M}_n} \exp\{-(K'\xi^2/2)(L_{m'})^{(1/2-\delta/2)_+}\}\nonumber\\
&&+ \sum_{m' \in \mathcal{M}_n}\frac{2\lambda_2\tilde\Gamma(m')}{n}\exp\{-(K'\xi^2/2)(L_{m'})^{(1/2-\delta/2)_+}\}
\label{I}.
\end{eqnarray}

\paragraph{\textbf{1) Case $0\leq \delta < 1/3$}} In that case, since
$\delta< (1/2-\delta/2)_+$, the choice $\xi^2=1$ ensures that
$\tilde \Gamma(m)\exp\{-(K'\xi^2/2)({L_m})^{(1/2-\delta/2)}\}$ is
bounded and thus the first term in (\ref{I}) is bounded by $C/n.$
Since $1\leq m\leq m_n$ with $m_n$ satisfying
(\ref{mn}),
$\sum_{m' \in \mathcal{M}_n}(\tilde\Gamma(m')/n)\exp\{-(K'/2)(L_{m'})^{(1/2-\delta/2)}\}
$ is bounded by $\tilde{\tilde C}/n,$
and hence $\sum_{m'\in {\mathcal M}_n}I(m^*)\leq C/n.$
Consequently, (\ref{but}) hold if we choose
$\mbox{pen}(m)=2a(1+2\xi^2)\lambda_1
({L_m})^{2\gamma+1-\delta}\exp\{2\mu\sigma^\delta(\pi
{L_m})^{\delta}\}/n. $
\paragraph{\textbf{2) Case $\delta= 1/3$}} According to the inequality
(\ref{eqmax}), $\xi^2$ is such that
$2\mu\sigma^\delta\pi^{\delta}(L_{m^*})^{\delta} -
(K'\xi^2/2)L_{m^*}^{\delta}= -2\mu\sigma^\delta (\pi
L_{m^*})^{\delta}$ that is
$\xi^2=(4\mu\sigma^\delta\pi^{\delta}\lambda_2)/(K_1\lambda_1).$
Arguing as for the case $0\leq \delta<1/3$, this
choice ensures that $\sum_{m'\in {\mathcal M}_n}I(m^*) \leq C/n$,
and consequently (\ref{but}) holds. The result follows for
$p(m,m')=2(1+2\xi^2)\lambda_1L_{m^*}^{2\gamma+1-\delta}\exp(2\mu\sigma^\delta(\pi
L_{m^*})^{\delta})/ n,$ and
$\mbox{pen}(m)=2a(1+2\xi^2)\lambda_1 {L_m}^{2\gamma+1-\delta}\exp(2\mu\sigma^\delta(\pi
{L_m})^{\delta})/ n.$

\paragraph{\textbf{3) Case $\delta> 1/3$}} If $\delta > (1/2-\delta/2)_+$,
according to (\ref{eqmax}) we choose
$\xi^2=\xi^2({L_m},L_{m'})$ such that
$2\mu\sigma^\delta\pi^{\delta}(L_{m^*})^{\delta} -
(K'\xi^2/2)L_{m^*}^{\omega}=
-2\mu\sigma^\delta\pi^{\delta}(L_{m^*})^{\delta}$ that is
$\xi^2=\xi^2(m,m')=(4\mu\sigma^\delta\pi^{\delta}\lambda_2)/(K_1\lambda_1)L_{m^*}^{\delta-\omega}.$
This choice ensures that $\sum_{m'\in {\mathcal M}_n}I(m^*) \leq
C/n$, and consequently (\ref{but}) holds
if
$p(m,m')=2(1+2\xi^2(m,m'))\lambda_1L_{m^*}^{2\gamma+1-\delta}\exp(2\mu\sigma^\delta(\pi
L_{m^*})^{\delta})/ n,$ associated to the penalty
$\mbox{pen}(m)=2a(1+2\xi^2({L_m},m))\lambda_1({L_m})^{2\gamma+1-\delta}\exp(2\mu\sigma^\delta(\pi
{L_m})^{\delta})/ n.$ \hfill $\Box$

\bigskip\noindent\textit{6.3 Technical Lemmas}
\newline

\noindent\textsc{Lemma 2}. \textit{
Let $\nu_n(t)$ be defined by (\ref{nu}),
$\Delta_1(m)$ be defined in Proposition~1. Under Assumptions
({\rm (A}$_1^{X,\varepsilon}${\rm )}
\begin{eqnarray}\label{b1}
\parallel\sum_{j\in \mathbb{Z}} |u^*_{\varphi_{m,j}}|^2\parallel_\infty \leq
\Delta_1(m),
\mbox{ and }
\sup_{ g\in \mathcal{S}_{s,r,b}(C_1)}\sum_{j\in
\mathbb{Z}} {\rm Var}[\nu_n(\varphi_{m,j})]\leq \Delta_1(m)/ n.
\end{eqnarray}} \newline

\paragraph{\textit{Proof of Lemma 2}}
Use the definition of $u^*_{\varphi_{m,j}}(z)$ to get that
\begin{eqnarray*}
\sum_{j \in \mathbb{Z}} \left\vert u^*_{\varphi_{m,j}}(z)\right\vert ^2=
\sum_{j \in \mathbb{Z}} \left\vert \int \exp\{ixz\} u_{\varphi_{m,j}}(x)
  dx\right\vert ^2
= \frac{L_m}{(2\pi)^2}\sum_{j \in \mathbb{Z}} \left\vert \int
\exp\{-ixz
L_m\}\exp\{ijx\}\frac{\varphi^*(x)}{f_\varepsilon^*(xL_m\sigma)}dx\right\vert
^2.
\end{eqnarray*}
By Parseval's Formula,
\begin{equation}\label{parseval}
\sum_{j \in \mathbb{Z}} \left\vert
u^*_{\varphi_{m,j}}(z)\right\vert ^2=(2\pi)^{-1} L_m\int
\left\vert
\frac{\varphi^*(x)}{f_\varepsilon^*(xL_m\sigma)}\right\vert^2dx=
\Delta_1(m),
\end{equation}
which entails that the first part of the bound (\ref{b1}) is
proved. The second part follows since
$\sum_{j\in \mathbb{Z}}\mbox{Var}[\nu_n(\varphi_{m,j})]
\leq   n^{-1}\int \sum_{j\in \mathbb{Z}}\left\vert
u^*_{\varphi_{m,j}}(z)\right\vert^2 h(z)dz.\;\; \Box
$ \newline

\noindent\textsc{Lemma 3}. \textit{
Let $\Delta_1(m)$  and $R(\mu,\delta,\sigma)$ be defined in Proposition~1
and in (\ref{lambda1}). Then
  under the assumption {\rm (A}$_2^{\varepsilon}${\rm )},
$\displaystyle \Delta_1(m)\leq
\frac 1{\pi \kappa_0^{2}R(\mu,\delta, \sigma)}(\pi
L_m)^{1-\delta}
(\sigma^2L_m^2\pi^2+1)^{\gamma}\exp\{2\mu\sigma^\delta\pi^{\delta}L_m^{\delta}\}.
$}\newline

\paragraph{\textit{Proof of Lemma 3}.}
Under the assumption {\rm (A}$_2^{\varepsilon}${\rm )}, $ \Delta_1(m)\leq (\pi
\kappa_0^{2})^{-1} (\sigma^2L_m^2\pi^2+1)^{\gamma}\int_0^{\pi L_m}
\exp\{2\mu\sigma^\delta u^{\delta}\} du.$ If $\delta=0$, by
convention $\mu=0$, and hence the integral in the previous bound
is less than $\pi L_m.$

Consider now the case $0<\delta \leq 1$. Easy calculations provide that
\begin{eqnarray*} \int_0^{\pi L_m}
e^{2\mu\sigma^\delta u^{\delta}} du &=& \int_0^{\pi L_m} \left(
  2\mu\sigma^\delta \delta
u^{\delta-1} e^{2\mu\sigma^\delta
  u^{\delta}}\right)\frac{du}{2\mu\sigma^\delta \delta
u^{\delta-1}} \leq  \frac{(\pi L_m)^{1-\delta}}{2\mu\sigma^\delta
  \delta}\left[e^{2\mu\sigma^\delta u^{\delta}}\right]_0^{\pi
L_m}\end{eqnarray*} and therefore $ \int_0^{\pi L_m}
\exp\{2\mu\sigma^\delta u^{\delta}\} du   \leq  [(\pi
L_m)^{1-\delta}/(2\mu\sigma^\delta\delta)]\exp(2\mu\sigma^\delta(\pi
L_m\sigma)^{\delta}). $

\noindent Now, if $\delta >1$, then by using that
$u^\delta=u^{\delta-1} u$, and consequently Lemma 3 follows from
 \begin{eqnarray*} \int_0^{\pi L_m} \exp\{2\mu\sigma^\delta
u^{\delta}\} du &\leq & \int_0^{\pi L_m} \exp\{2\mu\sigma^\delta(\pi
L_m)^{\delta-1} u\} du\leq\frac{(\pi
L_m)^{1-\delta}}{2\mu\sigma^\delta}\exp(2\mu\sigma^\delta(\pi
L_m)^{\delta}). \Box
\end{eqnarray*}

\noindent\textsc{Lemma 4}. \textit{
Let $\nu_n(t)$, $\Delta_1(m)$ and $\Delta_2(m, h)$ be
defined in (\ref{nu}), Proposition~1 and in (\ref{Delta2}). Then
under {\rm (A}$_1^{X,\varepsilon}${\rm )}
\begin{eqnarray*}
\sup_{ t \in B_{m,m'}(0,1) }\parallel u_t^*\parallel_\infty \leq
\sqrt{\Delta_1(m^*)}&&  \mathbb{E}[\sup_{t\in
B_{m,m'}(0,1)}|\nu_n(t)|]\leq \sqrt{\Delta_1(m^*)/n},\\
{\rm and } \sup_{t\in B_{m,m'}(0,1)} {\rm Var}(u^*_t(Z_1))&\leq &
\sqrt{\Delta_2(m^*,h)}/(2\pi).
\end{eqnarray*}}

\paragraph{\textit{Proof of Lemma 4}}
By combining Cauchy-Schwarz Inequality and (\ref{parseval}), the
square of the first term $\sup_{ t \in B_{m,m'}(0,1) }\parallel
u^*_t\parallel_\infty^2$ is bounded by $\sum_{j\in \mathbb{Z}}\int
\left\vert
    \varphi_{m^*,j}^*(u)/f_\varepsilon^*(\sigma u) \right\vert^2du=\Delta_1(m^*).$
Now, we have
$\mathbb{E}[\sup_{t\in B_{m,m'}(0,1)}|\nu_n(t)|]\leq
{\mathbb E}\left[(\sum_{j\in \mathbb{Z}}
(\nu_n(\varphi_{m^*,j}))^2)^{1/2}\right] \leq\left[\sum_{j\in \mathbb{Z}
  }\mbox{Var}(\nu_n(\varphi_{m^*,j}))\right]^{1/2},
$
which is bounded, by applying the second part of (\ref{b1}) in
Lemma 2, by $\sqrt{\Delta_1(m^*)/n}.$ Now write that
$
\sup_{ t\in B_{m,m'}(0,1)} \mbox{Var}(u^*_t(Z_1))\leq  \sup_{ t\in
  B_{m,m'}(0,1)} {\mathbb E}[\vert u^*_t(Z_1)\vert^2 ]
\leq  [\sum_{j,k\in \mathbb{Z}} \vert
Q_{j,k}(m^*)\vert^2]^{1/2},
$
with $ Q_{j,k}(m)={\mathbb E}[u^*_{\varphi_{m,j}}(Z_1) u^*_{\varphi_{m,k}}(-Z_1)]$
also given by
\begin{eqnarray*}
Q_{j,k}(m)=\frac{L_m}{(2\pi)^2}\iint \exp\{ijx-iky\}
\frac{\varphi^*(x)\varphi^*(y)}{f_\varepsilon^*(\sigma
L_mx)f_\varepsilon^*(\sigma L_my)}h^*(L_m(x-y))dxdy.
\end{eqnarray*}
Apply Parseval's Formula to get the result since
\begin{eqnarray*}
\sum_{j,k\in \mathbb{Z}} \vert Q_{j,k}(m)\vert^2&=&
\frac{L_m^2}{(2\pi)^2}\iint \left\vert
\frac{\varphi^*(x)\varphi^*(y)}{f_\varepsilon^*(\sigma
L_mx)f_\varepsilon^*(\sigma L_my)}h^*(L_m(x-y))
 \right\vert ^2dxdy. \Box \end{eqnarray*}

\bigskip\noindent\textsf{ACKNOWLEDGEMENTS}
We thank Cristina Butucea and Alexander Tsybakov for helpful discussions on
this subject.

\bigskip\noindent\textsf{REFERENCES}

\begin{description}
\item {\small A.R.\ Barron, L. Birg\'e \&  P. Massart (1999).
Risk bounds for model selection via penalization. \textit{
Probability Theory and Related  Fields}, 113, 301--413.}

\item {\small L. Birg\'e \&  Y. Rozenholc (2005). How many bins
must be put in a regular histogram. To appear in \textit{ ESAIM,
Probability and Statistics}.}

\item {\small C. Butucea (2004). Deconvolution of supersmooth
densities with   smooth noise. {\em  The Canadian Journal of
Statistics}, 32, 181--192.}

\item {\small  C. Butucea \&  A.B.\ Tsybakov (2004). Sharp
optimality and some effects of dominating bias in density
deconvolution, Preprint LPMA-898,
http://www.proba.jussieu.fr/mathdoc/preprints/index.html$\#$2004.}

\item {\small C. Butucea \&  C. Matias (2005). Minimax estimation
of the noise level and of the deconvolution density in a
semiparametric convolution model. \textit{ Bernoulli}, 11,  309-340.}

\item {\small R.J. Carroll, \&   P. Hall (1988). Optimal
rates of convergence for deconvolving a density. \textit{ Journal
of the American Statistical  Association}, 83, 1184--1186.}

\item {\small E.A. Cator (2001). Deconvolution with
  arbitrarily smooth kernels. \textit{ Statistics and Probability Letters}, 54,
  205--214.}

\item {\small L. Cavalier, Y. Golubev, O. Lepski \& A.B. Tsybakov~(2003).
Block  thresholding and sharp adaptive estimation in
severely ill-posed inverse problems.\textit{Theory Prob. Appl.},
48.}

\item {\small F. Comte and V. Genon-Catalot (2005) Penalized
projection estimator for volatility density. \textit{Preprint
MAP5} 2005-9,
http://www.math-info.univ-paris5.fr/map5/publis/titres05.html.}

 \item {\small A. Delaigle \&  I. Gijbels~(2004a). Practical
bandwidth selection in deconvolution kernel density estimation.
{\em Computational Statistics and Data Analysis}, 45, 249--267.}

\item {\small A. Delaigle \& I. Gijbels~(2004b). Bootstrap
bandwidth selection in kernel density estimation from a
contaminated sample. \textit{Annals of the Institute of
Statistical Mathematics}, 56, 19--47.}

\item {\small L. Devroye (1986). Non-Uniform Random Variable
Generation. Springer-Verlag, New-York.}

\item {\small L. Devroye (1989). Consistent deconvolution in
density estimation. {\em The Canadian Journal of Statistics}, 17,
235--239.}

\item {\small J. Fan (1991a). On the optimal rates of convergence
for nonparametric deconvolution problem. \textit{ The Annals of
Statistics}, 19,  1257--1272.}

\item {\small J. Fan (1991b). Asymptotic normality for
deconvolution kernel estimators, \textit{ Sankhya Series A}, 53,
97--110.}

\item {\small J. Fan \&  J.-Y. Koo (2002). Wavelet deconvolution.
\textit{ IEEE Transactions on Information Theory}, 48, 734--747.}

\item {\small C. Hesse (1999) Data-driven deconvolution. {\em
Journal of Nonparametric Statistics}, 10, 343--373.}

\item {\small  I.A. Ibragimov \&  R.Z. Hasminskii (1983).
Estimation of distribution density. \textit{ Journal of Soviet
Mathematics}, 21, 40--57.}

\item {\small  J.-Y Koo (1999). Logspline deconvolution in Besov
space. \textit{ Scandinavian Journal of Statistics}, 26, 73--86.}

\item {\small M.C.\ Liu \&  R.L.\ Taylor (1989). A consistent
nonparametric density estimator for the deconvolution problem.
\textit{ The Canadian Journal of Statistics}, 17, 427--438.}

\item {\small E.  Masry (1991). Multivariate probability density
deconvolution for stationary random processes. \textit{ IEEE
Transactions on Information Theory}, 37, 1105--1115.}

\item {\small A. Meister (2004). On the effect of misspecifying
the error density in a deconvolution problem. {\it The Canadian
Journal of Statistics}, 32(4), 439--449.}

\item {\small Y. Meyer (1990). Ondelettes et op\'erateurs, Tome I,
Hermann.}

\item {\small M. Pensky (2002). Density deconvolution based on
wavelets with bounded supports. \textit{ Statistics and
Probability Letters}, 56, 261--269.}

\item {\small M. Pensky \&  B. Vidakovic (1999). Adaptive wavelet
estimator for nonparametric density deconvolution. \textit{ The
Annals of Statistics}, 27, 2033--2053.}

\item {\small L. Stefanski \&  R.J. Carroll (1990). Deconvoluting
kernel density estimators. {\em Statistics}, 21, 1696-184.}

\item {\small M. Talagrand (1996). New concentration inequalities
in product spaces. \textit{ Inventiones Mathematicae}, 126,
505--563.}

\item {\small A.B. Tsybakov (2000). On the best rate of adaptive
estimation in some inverse problems. \textit{ Comptes-Rendus de
l'Acad\'emie des Sciences, Paris, S\'erie I Math\'ematiques}, 330,
835-840.}

\item {\small C.H. Zhang (1990). Fourier methods for estimating
mixing densities and distributions. \textit{ The Annals of
Statistics}, 18, 806-831.}

\end{description}


\hrule
\newpage

\hfill{Fabienne COMTE: \texttt{fabienne.comte@univ-paris5.fr}}%
\begin{flushright}
\vglue -2.3mm
{\it Universit\'e Ren\'e Descartes-Paris 5}\\
{\it 45, rue des Saint-P\`eres}\\
{\it 75270 Paris Cedex 06 FRANCE}\\
\bigskip
Yves ROZENHOLC: {\tt yves.rozenholc@math-info.univ-paris5.fr}\\
{\it Laboratoire MAP5, UMR 5145}\\
{\it Universit\'e Ren\'e Descartes-Paris 5}\\
{\it 45, rue des Saint-P\`eres}\\
{\it 75270 Paris Cedex 06 FRANCE}\\
\bigskip
Marie-Luce TAUPIN: {\tt marie-luce.taupin@math.u-psud.fr}\\
{\it D\'epartement de math\'ematiques, UMRC 8628\\
Universit\'e Paris-Sud\\
{\it 91405 Orsay FRANCE}}
\end{flushright}

\end{document}